# READING POLICIES FOR JOINS: AN ASYMPTOTIC ANALYSIS

By Ralph P. Russo and Nariankadu D. Shyamalkumar[1]

*University of Iowa*

Suppose that $m_n$ observations are made from the distribution $\mathbf{R}$ and $n - m_n$ from the distribution $\mathbf{S}$. Associate with each pair, $x$ from $\mathbf{R}$ and $y$ from $\mathbf{S}$, a nonnegative score $\phi(x, y)$. An *optimal* reading policy is one that yields a sequence $m_n$ that maximizes $\mathbb{E}(M(n))$, the expected sum of the $(n - m_n)m_n$ observed scores, *uniformly* in $n$. The *alternating* policy, which switches between the two sources, is the optimal nonadaptive policy. In contrast, the *greedy* policy, which chooses its source to maximize the expected gain on the next step, is shown to be *the* optimal policy. Asymptotics are provided for the case where the $\mathbf{R}$ and $\mathbf{S}$ distributions are discrete and $\phi(x, y) = 1$ or 0 according as $x = y$ or not (i.e., the observations match). Specifically, an invariance result is proved which guarantees that for a wide class of policies, including the *alternating* and the *greedy*, the variable $M(n)$ obeys the same CLT and LIL. A more delicate analysis of the sequence $\mathbb{E}(M(n))$ and the sample paths of $M(n)$, for both *alternating* and *greedy*, reveals the slender sense in which the latter policy is asymptotically superior to the former, as well as a sense of equivalence of the two and robustness of the former.

**1. Introduction.** Suppose that samples of size $m_n$ and $n - m_n$ are drawn from tables $\mathbf{R}$ and $\mathbf{S}$ in a database, table $\mathbf{R}$ containing information (age, interests, education level, etc.) on a group of single males, and table $\mathbf{S}$ the same information on a group of single females. Associate with each pair of records, $x$ from $\mathbf{R}$ and $y$ from $\mathbf{S}$, a nonnegative score $\phi(x, y)$ whose value depends on how closely the two records agree. A male and female of similar age, with common interests and education level, would have a high score (a value near 1 on a $[0, 1]$ scale, e.g.). The goal is to choose $m_n$ to maximize $\mathbb{E}(M(n))$, where $M(n)$ is the sum of the $m_n(n - m_n)$ scores generated by the $n$ records that have been read. In this way, the expected overall interest

Received February 2006; revised July 2006.

[1]Supported in part by the University of Iowa Old Gold Fellowship.

level between the two groups (after $n$ reads) is maximized. Alternatively, **R** may contain information on a group of buyers in a marketplace (specifically, which items each seeks to buy) and **S** information on a group of sellers (which items each seeks to sell), the goal then being to maximize the level of commerce between the groups.

*The alternating and myopic policies.* Suppose that observations are made sequentially and without replacement from each of two sources (populations) **R** and **S**. An algorithm that sequentially chooses the source for each observation is referred to as a reading policy. An optimal reading policy (if existent) is one that maximizes $\mathbb{E}(M(n))$ uniformly in $n$. Two reading policies of interest are the *alternating*, which alternately samples from **R** and **S**, and the *myopic* (or *greedy*), which on each step chooses the source that maximizes the expected gain $\mathbb{E}(M(n)) - \mathbb{E}(M(n-1))$ for that step. The alternating policy is interesting because it is easy to implement and requires no knowledge of **R** or **S**. Moreover, this policy is optimal in a restricted sense (see below). Any policy with a fixed sampling order for which the **R** sample size is always within one of the **S** sample size is considered an alternating policy, as all such policies produce the same expected total score at all steps. In contrast to the alternating policy, the greedy policy requires a complete knowledge of **R** and **S**. It is a short term strategy that optimizes the expected gain on the next step, with no explicit regard to future gains. Note that there may be more than one greedy policy, as occasionally the greedy criterion may be ambivalent between **R** and **S**.

In the case of the equijoin, the records $x$ from **R** and $y$ from **S** can be categorized by positive integer values $r(x)$ and $s(y)$ with $\phi(x, y) = 1$ or $0$ accordingly as $r(x) = s(y)$ or not. When $\phi(x, y) = 1$, we say that records $x$ and $y$ *match*. Optimality in the case of the *equijoin* was studied in [16]. When **R** or **S** is finite, it was shown that an optimal policy need not exist, that the alternating policies are optimal among the restricted class of nonadaptive policies (those that ignore the information obtained from the samples), and that any greedy policy dominates (and in most cases is strictly better than) any alternating policy. That *alternating* is the optimal *nonadaptive* (**R** and **S** both infinite or not and $\phi$ arbitrary) and is easy to show, so is stated here without proof.

When **R** and **S** are infinite, the problem reduces to i.i.d. sampling from those distributions. In this case it is shown in [16] that *alternating* is again optimal among the *nonadaptives*, and that *greedy* is optimal among *all* reading policies. In the next section we provide a simpler proof of a much stronger result; namely, that *greedy* is the optimal policy under the so-called total expected discounted reward criterion, for any decreasing discount sequence. From this it follows that *greedy* is the optimal policy in the $\phi$ arbitrary case. This case includes an interesting class of score functions which satisfy $\phi(x, y) = 1$ or $0$ (like the equijoin), but (unlike the equijoin) allows $\phi(x_1, y_1) = \phi(x_2, y_1) = \phi(x_2, y_2) \neq \phi(x_1, y_2)$. Such is the case when all observations are points on a space and $\phi(x, y) = 1$ or $0$ accordingly as $x$ from **R** and $y$ from **S** are within a prescribed distance $t$ of each other.

Our interest in the *alternating* and *greedy* policies stems from the above optimality properties. Our main focus is on the asymptotic properties of the alternating and greedy policies in the i.i.d. (infinite populations) case. To simplify the presentation, we confine our attention to the *equijoin*. A preview of the kinds of results we seek is provided in the following example. Interestingly, even in this simple scenario, the analysis is not trivial.

*An illustrative example.* On each step of a coin tossing experiment suppose that one may choose either of two coins to toss: one of them fair and the other two-headed. Let $M_{\mathbf{A}}(n)$ denote the numbers of matches formed by the policy that alternates between the coins (starting with the fair) after a total of $n$ tosses ($n$th epoch) have been made. We have that $M_{\mathbf{A}}(n) = (n/2) \operatorname{Bin}(n/2; 1/2)$ for $n$ even and $M_{\mathbf{A}}(n) = [(n-1)/2] \operatorname{Bin}((n+1)/2; 1/2)$ for $n$ odd; $\operatorname{Bin}(m; p)$ being a binomial random variable with parameters $m \geq 1$ and $p \in [0, 1]$. In particular, this implies that at the $n$th epoch the expected numbers of matches equals $n^2/8$ or $(n^2 - 1)/8$ accordingly as $n$ is even or odd and that

$$\left( \frac{M_{\mathbf{A}}(n) - n^2/8}{n^{3/2}} \right) \xrightarrow{\mathrm{d}} \mathrm{N}\left( 0, \frac{1}{32} \right).$$

It can be easily checked that the following is a member of the class of greedy policies and therefore optimal: *toss the fair coin until heads is obtained, toss the two-headed coin twice, return to the fair coin and repeat the cycle*. We denote the number of matches at the $n$th epoch using this greedy policy by $M_{\mathbf{G}}(n)$.

The derivation of a closed form expression for $\mathbb{E}(M_{\mathbf{G}}(n))$ is a bit more involved than it was for $\mathbb{E}(M_{\mathbf{A}}(n))$. A method outlined in Section 4 yields

(1.1)
$$\begin{aligned}
&\mathbb{E}(M_{\mathbf{G}}(n)) \\
&= \frac{n^2}{8} + \frac{n}{16} - \frac{7}{64} \\
&\quad + \frac{3\sin((n-1)\beta) + 16\sin((n-3)\beta) - 9\sqrt{7}\cos((n-1)\beta)}{2^{(n+11)/2}\sqrt{7}},
\end{aligned}$$

where $\beta := \pi - \arctan(\sqrt{7})$. In particular, this implies that

(1.2)
$$\left( \frac{1}{n} \right) \mathbb{E}(M_{\mathbf{G}}(n) - M_{\mathbf{A}}(n)) \to \frac{1}{16},$$

which in turn implies that $\mathbb{E}(M_{\mathbf{A}}(n)) < \mathbb{E}(M_{\mathbf{G}}(n)) < \mathbb{E}(M_{\mathbf{A}}(n+1))$ for $n \geq 3$. Thus, there exists a rather tight link between the expectations of the two processes.

An approach to understanding $M_{\mathbf{G}}(n)$ [and not just $\mathbb{E}(M_{\mathbf{G}}(n))$] uses an embedded renewal counting process $\{N(n)\}_{n \geq 1}$, with a renewal occurring upon the observance of a tail, and an inter-arrival variable $3Z + 1$, where $Z$ is a unit mean geometric random variable. The relation between $M_{\mathbf{G}}(n)$ and $N(n)$ depends on the state occupied at the $n$th epoch; the states being a tail, heads with fair coin, first heads (i.e., not preceded by another) with 2-headed coin and second heads (i.e., preceded by another) with 2-headed coin. For example, $M_{\mathbf{G}}(n) = (2/9)(n - N(n) + 2)(n - N(n) - 1)$ when the process has just observed a tail. We note that the state process is a doubly stochastic Markov chain.

The above with the approximation $M_{\mathbf{G}}(n) \approx (2/9)(n - N(n))^2$ and the CLT for renewal counting processes (see [12], page 62) implies that $M_{\mathbf{G}}(n)$ has the same weak limit as $M_{\mathbf{A}}(n)$. Law of the iterated logarithm results for $M_{\mathbf{G}}(\cdot)$ and $M_{\mathbf{A}}(\cdot)$ are likewise easily obtainable and again coincide. Also, the exact expressions relating $N(n)$ and $M_{\mathbf{G}}(n)$ along with the geometric rate of convergence to stationarity of the state process Markov chain and the expectation of the residual lifetime (or *overshoot*) of the renewal process yields (1.2).

Another phenomenon that we find interesting is the following: when both policies are driven by the same sequence of fair coin tosses, *alternating* beats (produces more matches than) *greedy* on infinitely many epochs, with probability one—the optimality of the *greedy* notwithstanding. The reason for this is that the coin sequence is obliged to transition (in two steps) from the state ($k$ heads, $k$ tails) to ($k$ heads, $k + 2$ tails) infinitely often. And when it does, we observe that *alternating* produces $k$ more matches than *greedy* upon completion of the $4k + 2$nd step. It can similarly be argued that *greedy* beats *alternating* infinitely often with probability one.

1.1. *Overview of results.* Using a method from dynamic programming, we prove in Section 2 that the greedy policy is optimal under the total expected reward criterion for any finite horizon. This result extends our result in [16] (that greedy is optimal in the i.i.d. equijoin case) to general score functions $\phi$. For our asymptotic analysis, rather than exploit a renewal structure (as in the example above), we instead take advantage of an embedded martingale structure in order to produce a broader range of results.

A key to the weak and strong limiting behavior of $M_{\mathbf{G}}(n)$ is an invariance result (proved in Section 3) which says that the asymptotic behavior of $M(n)$ under any policy is governed by the variable $R(n)$, the number of records read from $\mathbf{R}$ from among the first $n$ records read. We prove a central limit theorem and law of the iterated logarithm for $R_{\mathbf{G}}(n)$, which yields (by invariance) a common CLT and LIL for $M_{\mathbf{A}}(n)$ and $M_{\mathbf{G}}(n)$. Thus, an observer working with perfect knowledge of the distributions of $\mathbf{R}$ and $\mathbf{S}$ can not do much better (produce more matches) using the greedy policy than his counterpart who uses the alternating policy and is ignorant of those distributions.

In Section 4 we take up the mathematical question of how much better is *greedy* than *alternating*. Specifically, we find an expression for $\mathbb{E}(M_{\mathbf{G}}(n))$ that is similar to (1.1), but which contains a low order (linear) error term. As in the illustrative example, this yields

$$\lim_{n \to \infty} \frac{\mathbb{E}(M_{\mathbf{A}}(n) - M_{\mathbf{G}}(n))}{n} > 0$$

and for a finite constant $k$ computable from the distributions of $\mathbf{R}$ and $\mathbf{S}$,

$$\mathbb{E}(M_{\mathbf{A}}(n)) < \mathbb{E}(M_{\mathbf{G}}(n)) < \mathbb{E}(M_{\mathbf{A}}(n+k)) \qquad \text{for all large } n.$$

The former statement uncovers a measure by which *greedy* is asymptotically superior to *alternating*, while the latter reveals how tightly connected the two processes are, in terms of their expectations. We next identify the weak limit of $(M_{\mathbf{G}}(n) - M_{\mathbf{A}}(n))/n^{5/4}$ as a scale mixture of normals, showing that $M_{\mathbf{G}}(n)$ and $M_{\mathbf{A}}(n)$ differ by a higher order (than linear) stochastic term which is symmetric about zero. Finally, we present a crude LIL type result for $M_{\mathbf{G}}(n) - M_{\mathbf{A}}(n)$ which shows that although $\mathbb{E}(M_{\mathbf{A}}(n))$ and $\mathbb{E}(M_{\mathbf{G}}(n))$ are tightly linked via the above inequality, it takes an arbitrarily large number of epochs (infinitely often, with probability one) for the sample path of one process to catch up with that of the other.

1.2. *Notation.* All vectors carry a tilde. For two vectors, $\tilde{x}$ and $\tilde{y}$, with the same dimensions, $\tilde{x} \cdot \tilde{y}$ will denote their inner product. Almost sure convergence, convergence in probability and weak convergence will be denoted by $\xrightarrow{\text{a.s.}}$, $\xrightarrow{\text{P}}$ and $\xrightarrow{\text{d}}$, respectively. We denote the iterated logarithm by $\log_2$, that is, $\log_2(n) = \log(\log(n))$.

2. **Preliminaries.** We consider two sources $\mathbf{R}$ and $\mathbf{S}$ with both containing infinitely many records. A record from either source, $\mathbf{R}$ or $\mathbf{S}$, carries a single positive integer valued label. The probability that a record from the $\mathbf{R}$ source (resp., the $\mathbf{S}$ source) carries the $i$th label is $r_i$ (resp., $s_i$). The probability vectors $(r_1, r_2, \ldots)$ and $(s_1, s_2, \ldots)$ are denoted by $\tilde{r}$ and $\tilde{s}$, respectively. The inner product of $\tilde{r}$ with $\tilde{s}$ is denoted by $\mu$, that is, $\mu = \tilde{r} \cdot \tilde{s}$. We shall assume $\mu$ to be positive, as otherwise there will be no common label between the two sources. The labels on the $n$th records read from the $\mathbf{R}$ and $\mathbf{S}$ sources are denoted by $L_{\mathbf{R}}(n)$ and $L_{\mathbf{S}}(n)$, respectively. The above implies that $\{L_{\mathbf{R}}(n)\}_{n \geq 1}$ and $\{L_{\mathbf{S}}(n)\}_{n \geq 1}$ are sequences of independent and identically distributed random variables with

$$\Pr(L_{\mathbf{R}}(1) = i) = r_i \quad \text{and} \quad \Pr(L_{\mathbf{S}}(1) = i) = s_i, \qquad i = 1, 2, \ldots.$$

Associated with the sequences $\{L_{\mathbf{R}}(n)\}_{n\geq 1}$ and $\{L_{\mathbf{S}}(n)\}_{n\geq 1}$ are the discrete time vector counting processes $\{\tilde{N}_{\mathbf{R}}(n)\}_{n\geq 1}$ and $\{\tilde{N}_{\mathbf{S}}(n)\}_{n\geq 1}$; the first is defined by

$$\tilde{N}_{\mathbf{R}}(n) = (N_{\mathbf{R}}(n,1), N_{\mathbf{R}}(n,2), \ldots),$$

$$\text{with } N_{\mathbf{R}}(n,i) = \sum_{j=1}^{n} I_{\{L_{\mathbf{R}}(j)=i\}}, \ i, n \geq 1,$$

and the second is defined analogously.

2.1. *Reading policies.* A reading policy is a zero–one valued stochastic process

$$C(n) = \begin{cases} 1, & \text{if the } n\text{th selection is from } \mathbf{R}, \\ 0, & \text{if the } n\text{th selection is from } \mathbf{S}, \end{cases} \quad n = 1, 2, \ldots.$$

Associated with each reading policy are two counting processes $\{R(n)\}_{n\geq 1}$ and $\{S(n)\}_{n\geq 1}$ defined by

$$R(n) := \sum_{j=1}^{n} C(j) \quad \text{and} \quad S(n) := n - R(n), \qquad n = 1, 2, \ldots.$$

These processes keep track of the number of records read from $\mathbf{R}$ and $\mathbf{S}$, respectively. We shall refer to $R(n)/n$ as the *selection ratio*. Also associated with a reading policy is a nondecreasing process $\{M(n)\}_{n\geq 1}$ which counts the number of matches, generated by the first $n$ records:

$$M(n) = \tilde{N}_{\mathbf{R}}(R(n)) \cdot \tilde{N}_{\mathbf{S}}(S(n)), \qquad n = 1, 2, \ldots.$$

Observe that all of the processes $\{M(n)\}_{n\geq 1}$, $\{R(n)\}_{n\geq 1}$ and $\{S(n)\}_{n\geq 1}$ depend on the reading policy even though the notation does not make it explicit.

The filtration $\{\mathcal{F}_n\}_{n\geq 0}$ for a given reading policy is defined by

$$\mathcal{F}_n := \mathcal{F}_0 \vee \sigma \langle L_{\mathbf{R}}(1), \ldots, L_{\mathbf{R}}(R(n)); L_{\mathbf{S}}(1), \ldots, L_{\mathbf{S}}(S(n)) \rangle, \qquad n = 1, 2, \ldots,$$

with $\mathcal{F}_0$ containing all the information needed for randomization and independent of $\{L_{\mathbf{R}}(n)\}_{n\geq 1}$ and $\{L_{\mathbf{S}}(n)\}_{n\geq 1}$. All reading policies are required to be predictable with respect to the above filtration—otherwise they would not be *implementable*.

DEFINITION 2.1. An alternating policy is a $\mathcal{F}_0$ measurable reading policy for which

$$(2.1) \qquad\qquad R(2n) = n, \qquad n = 1, 2, \ldots.$$

In words, an alternating policy is one which does not use any information from the records, and under which at any step the numbers of records read from the two sources are within one of each other. There exists an infinite number of alternating policies. One of the simplest alternating policies is defined by $C(n) = n \bmod 2$. In fact, in the arguments we tacitly assume for convenience that we are working with this version. From the point of view of implementation though, one may prefer the alternating policy given by $C(n) = I_{\{n \bmod 4 < 2\}}$ as it, leaving apart the first record, reads two records at a time from the chosen source.

Toward defining greedy policies, we observe that

$$
\begin{aligned}
\mathbb{E}(M(n+1) - M(n)|\mathcal{F}_n) = {} & \mathbb{E}(N_{\mathbf{S}}[S(n), L_{\mathbf{R}}(R(n)+1)]|\mathcal{F}_n)C(n+1) \\
& + \mathbb{E}(N_{\mathbf{R}}[R(n), L_{\mathbf{S}}(S(n)+1)]|\mathcal{F}_n)(1 - C(n+1)).
\end{aligned}
$$

Hence, any reading policy $C(\cdot)$ maximizing the above conditional expectation should satisfy, for $n \geq 1$,

$$
(2.2) \quad C(n+1) = \begin{cases}
1, & \text{if } \mathbb{E}(N_{\mathbf{S}}[S(n), L_{\mathbf{R}}(R(n)+1)]|\mathcal{F}_n) \\
& \quad > \mathbb{E}(N_{\mathbf{R}}[R(n), L_{\mathbf{S}}(S(n)+1)]|\mathcal{F}_n), \\
0, & \text{if } \mathbb{E}(N_{\mathbf{S}}[S(n), L_{\mathbf{R}}(R(n)+1)]|\mathcal{F}_n) \\
& \quad < \mathbb{E}(N_{\mathbf{R}}[R(n), L_{\mathbf{S}}(S(n)+1)]|\mathcal{F}_n),
\end{cases}
$$

with no requirement on epochs where

$$
(2.3) \quad \mathbb{E}(N_{\mathbf{S}}[S(n), L_{\mathbf{R}}(R(n)+1)]|\mathcal{F}_n) = \mathbb{E}(N_{\mathbf{R}}[R(n), L_{\mathbf{S}}(S(n)+1)]|\mathcal{F}_n).
$$

As $\{L_{\mathbf{R}}(n)\}_{n \geq 1}$ and $\{L_{\mathbf{S}}(n)\}_{n \geq 1}$ are sequences of i.i.d. random variables, we have

$$
\mathbb{E}(N_{\mathbf{S}}[S(n), L_{\mathbf{R}}(R(n)+1)]|\mathcal{F}_n) = \tilde{N}_{\mathbf{S}}(S(n)) \cdot \tilde{r}, \qquad n = 1, 2, \ldots,
$$

and an analogous relation for $\tilde{N}_{\mathbf{R}} R(n) \cdot \tilde{s}$.

Our analysis depends on the observation that $\tilde{N}_{\mathbf{S}} S(n) \cdot \tilde{r}$ and $\tilde{N}_{\mathbf{R}} R(n) \cdot \tilde{s}$ are both partial sums of i.i.d. observations. To make this explicit, we define

$$
X_{\mathbf{R}}(n) := s_{L_{\mathbf{R}}(n)} \quad \text{and} \quad X_{\mathbf{S}}(n) := r_{L_{\mathbf{S}}(n)}, \qquad n = 1, 2, \ldots.
$$

The two sequences $\{X_{\mathbf{R}}(n)\}_{n \geq 1}$ and $\{X_{\mathbf{S}}(n)\}_{n \geq 1}$ are sequences of i.i.d. random variables with common mean $\mu$ and variances $\sigma_{\mathbf{R}}^2$ and $\sigma_{\mathbf{S}}^2$, respectively. We shall denote their partial sums by $\Gamma_{\mathbf{R}}[\cdot]$ and $\Gamma_{\mathbf{S}}[\cdot]$, that is,

$$
\Gamma_{\mathbf{R}}[n] = \sum_{j=1}^{n} X_{\mathbf{R}}(j) \quad \text{and} \quad \Gamma_{\mathbf{S}}[n] = \sum_{j=1}^{n} X_{\mathbf{S}}(j), \qquad n = 1, 2, \ldots.
$$

This leads to the relations

$$
\tilde{N}_{\mathbf{S}}(S(n)) \cdot \tilde{r} = \Gamma_{\mathbf{S}}[S(n)] \quad \text{and} \quad \tilde{N}_{\mathbf{R}}(R(n)) \cdot \tilde{s} = \Gamma_{\mathbf{R}}[R(n)], \qquad n = 1, 2, \ldots.
$$

Combining the above with (2.2) leads to the following definition.

DEFINITION 2.2. A reading policy $C(\cdot)$ is called a *greedy* policy if it satisfies

$$(2.4) \quad C(n+1) = \begin{cases} 1, & \text{if } \Gamma_{\mathbf{S}}[S(n)] > \Gamma_{\mathbf{R}}[R(n)], \\ 0, & \text{if } \Gamma_{\mathbf{S}}[S(n)] < \Gamma_{\mathbf{R}}[R(n)], \end{cases} \qquad n = 1, 2, \ldots.$$

Henceforth, all quantities with a subscript of $\mathbf{G}$ will pertain to a greedy policy and those with a subscript of $\mathbf{A}$ to an alternating policy. An important consequence of the definition of a greedy policy is that

$$(2.5) \quad |\Gamma_{\mathbf{R}}[R_{\mathbf{G}}(n)] - \Gamma_{\mathbf{S}}[S_{\mathbf{G}}(n)]| \leq \gamma, \qquad n = 1, 2, \ldots,$$

where $\gamma := \max_{1 \leq i \leq \infty} r_i \vee \max_{1 \leq i \leq \infty} s_i$.

The case $\sigma_{\mathbf{R}} + \sigma_{\mathbf{S}} = 0$ is equivalent to having $\tilde{r}$ and $\tilde{s}$ as uniform distributions with identical finite supports. And it is easily checked that identical uniform distributions make the set of all greedy policies coincide with the set of all alternating policies. Hence, in the following we will assume $\sigma_{\mathbf{R}} + \sigma_{\mathbf{S}} > 0$.

2.2. *Optimality of the greedy.* Here we show using dynamic programming that the greedy policy maximizes the expected number of matches at all epochs in the case of infinite populations. A key observation to showing this is that the incremental gain of matches from the $(n+1)$st record can be written in terms of $\Gamma_{\mathbf{R}}[R(n)]$ and $\Gamma_{\mathbf{S}}[S(n)]$ as

$$\mathbb{E}(M(n+1) - M(n)|\mathcal{F}_n) = \Gamma_{\mathbf{S}}[S(n)]C(n+1) + \Gamma_{\mathbf{R}}[R(n)](1 - C(n+1)).$$

This suggests a rather compact Markov Decision Problem (MDP) formulation—at the $n$th epoch, the state is defined as $(\Gamma_{\mathbf{R}}[R(n)], \Gamma_{\mathbf{S}}[S(n)])$ and the action of choosing the next record from the $\mathbf{R}$ source results in a reward of $\Gamma_{\mathbf{S}}[S(n)]$ with $(\Gamma_{\mathbf{R}}[R(n)] + X_{\mathbf{R}}(R(n) + 1], \Gamma_{\mathbf{S}}[S(n)])$ as the new state (when the next record is chosen from $\mathbf{S}$ the reward and the new state are analogously defined). That this compact representation fails in the case of finite population(s) is easily demonstrated; see [16].

Abstractly, following the system of specifying a MDP as given in [11], consider the MDP with decision epochs $\{1, \ldots, N\}$ for some $N \geq 1$, state space $[0, \infty) \times [0, \infty)$, action set (invariant to the current state) $\{0, 1\}$, with time homogeneous expected rewards

$$r(\tilde{xi}, a) = \xi_{1+a} \qquad \text{for } a = 0, 1; \tilde{\xi} = (\xi_1, \xi_2) \in [0, \infty) \times [0, \infty)$$

and time homogeneous transition probabilities

$$p(\tilde{\xi}'|\xi; a) = \begin{cases} p_2(\xi_2' - \xi_2), & a = 0 \text{ and } \xi_1' = \xi_1, \\ p_1(\xi_1' - \xi_1), & a = 1 \text{ and } \xi_2' = \xi_2, \\ 0, & \text{otherwise}, \end{cases}$$

where $p_1(\cdot)$ and $p_2(\cdot)$ are probability densities (with respect to some $\sigma$-finite measure $\lambda$) on $[0, \infty)$ with a common mean, say, $\theta$. In terms of our original problem, action 1 (resp. 0) corresponds to picking the next record from **R** (resp. **S**), $p_1(\cdot)$ [resp., $p_2(\cdot)$] corresponds to the mass function of $X_\mathbf{R}$ (resp. $X_\mathbf{S}$) and the reward is the expected increment in the number of matches from the next record.

The above MDP, while reminiscent of a two-armed bandit (see, e.g., [11] and [13]), is not quite so. Considering the renewal processes, with inter-arrival distributions $p_1(\cdot)$ and $p_2(\cdot)$, as the states of two projects, the expected reward generated by choosing a project is equal to the state of the other. This dependence of the reward on the state of the other (idle) project fails one of the requirements of the bandit problem; see [11]. Nevertheless, it fits the formalism of the (two machine) tax problem of [18] (ongoing bandits of [4]) of where the reward structure is in a sense the reverse of the bandit problem.

In the tax problem, at any epoch, one of $K$ machines can be operated with idle machines generating a cost (the *tax*). The goal is to schedule the machines in order to minimize, for example, the expected total discounted stream of costs. Interestingly, from the point of view of a search for the optimal strategy, the tax problem (with infinite horizon and discounted rewards) and the bandit problem are equivalent; see [18]. As shown in [1], such an equivalence holds even while allowing all machines, active or inactive, to generate either a reward or a cost (negative reward) with the goal of maximizing the total discounted stream of rewards. Also, our MDP can be seen to be a particular case of the generalized bandit problem of [10]. While [1, 10] and [18] look at an infinite horizon discounted reward criterion, our interest is in the finite horizon analysis of our MDP. Below we show by a simple inductive argument that the greedy (myopic) policy is optimal in the finite horizon case under the total expected reward criterion. For an involved proof using the interchange argument; see [16]. Also, it is not hard to construct a simple qualitative argument along the lines of the proof of the Gittin's index theorem of [19].

It should be no surprise, given the time homogeneity and two point action set, that optimal deterministic Markov policies exist; see, for example, Theorem 4.4.2 of [11]. Moreover, it is easily argued from the reward structure that an optimal deterministic Markov policy which is a function of $\xi_1 - \xi_2$ exists. Below we additionally show that this policy is given by $I_{(-\infty,0)}(\xi_1 - \xi_2)$, the greedy policy.

THEOREM 2.1. *The greedy policy is optimal under the total expected reward criterion for any finite horizon.*

REMARK 2.1. That the greedy policy maximizes the expected number of matches at all epochs implies that it also maximizes the total expected discounted incremental matches for all nonincreasing discount sequences.

REMARK 2.2. Note that the above theorem also implies that the greedy maximizes

$$\mathbb{E}\left(\sum_{i=1}^{R(n)} \sum_{j=1}^{S(n)} \phi(L_{\mathbf{R}}(i), L_{\mathbf{S}}(j))\right),$$

where $\phi(\cdot, \cdot)$ is nonnegative. This implies that optimality of the greedy extends beyond *equijoins*. Moreover, the proof allows $L_{\mathbf{R}}$ and $L_{\mathbf{S}}$ to be random elements on any space.

REMARK 2.3. The above theorem does not extend directly beyond two sources. This is reminiscent of the bandit problem with different discount factors for each bandit—Gittin's index exists in the case of the two armed bandit, but not beyond. For details, we refer to [4].

PROOF OF THEOREM 2.1. The proof is by induction. Let $V_n(\xi)$ denote the maximum total expected reward for the $n$-epoch ($n$ more records to pick) problem at state $\xi$. Clearly, $V_1(\xi) = \max(\xi_1, \xi_2)$, which is attained by the greedy policy. Assume without loss of generality that $\xi_1 \geq \xi_2$. Now, since

$$\begin{aligned}
\xi_1 + \int V_1((\xi_1, \xi_2 + \zeta)) \, dp_2(\zeta) &\geq \xi_1 + \int (\xi_2 + \zeta) \, dp_2(\zeta) \\
&= \xi_1 + \xi_2 + \theta \\
&= \xi_2 + \int (\xi_1 + \zeta) \, dp_1(\zeta) \\
&= \xi_2 + \int V_1((\xi_1 + \zeta, \xi_2)) \, dp_1(\zeta),
\end{aligned}$$

we have the greedy policy is optimal for the 2-epoch problem too. Now assume that the greedy policy attains $V_i(\xi)$ for $i = 1, 2, \ldots, (n-1)$ for all $\xi \in [0, \infty) \times [0, \infty)$. That $V_n(\xi)$ is also attained by a greedy follows from

$$\begin{aligned}
&\xi_1 + \int V_{n-1}((\xi_1, \xi_2 + \zeta)) \, dp_2(\zeta) \\
&\geq \xi_1 + \int \left[ \int \xi_2 + \zeta + V_{n-2}((\xi_1 + \vartheta, \xi_2 + \zeta)) \, dp_1(\vartheta) \right] dp_2(\zeta) \\
&= \xi_1 + \xi_2 + \theta + \int \int V_{n-2}((\xi_1 + \vartheta, \xi_2 + \zeta)) \, dp_1(\vartheta) \, dp_2(\zeta) \\
&= \xi_2 + \int \left[ \int \xi_1 + \vartheta + V_{n-2}((\xi_1 + \vartheta, \xi_2 + \zeta)) \, dp_2(\zeta) \right] dp_1(\vartheta) \\
&= \xi_2 + \int V_{n-1}((\xi_1 + \vartheta, \xi_2)) \, dp_1(\vartheta).
\end{aligned}$$

Hence, the proof. $\quad\square$

**3. Basic weak and strong limit theorems.** The heuristics in the introduction suggest (and the results of this section confirm) that for a policy to be competitive, its selection ratio must converge to $1/2$. However, in some applications the observer may not have control of the sampling order, or may find it cost effective to sample unevenly from the two sources. Thus, the case where $R(n)/n \xrightarrow{\text{a.s.}} \alpha \in (0,1)$ is of interest. In this section we show that the number of matches, when suitably centered and scaled, can be strongly approximated by a standard Wiener process. From this approximation it is easy to obtain a CLT and LIL for the number of matches.

THEOREM 3.1. *Consider a reading policy that satisfies $R(n)/n \xrightarrow{\text{a.s.}} \alpha \in (0,1)$ and the associated process $Z(\cdot)$ defined by*

$$Z(t) := \begin{cases} 0, & \text{for } 0 \le t < V_1, \\ \dfrac{M(n)}{n} - n\alpha(1-\alpha)\mu, & \text{for } V_n \le t < V_{n+1} \text{ and } n \ge 1, \end{cases}$$

*where $V_n := (1-\alpha)^2 \sigma_R^2 R(n) + \alpha^2 \sigma_S^2 S(n)$ for $n \ge 1$. If for some $\beta \ge 0$,*

$$(3.1) \qquad \sqrt{\frac{n}{[\log_2(n)]^{(1-\beta)}}} \left( \frac{R(n)}{n} - \alpha \right) \xrightarrow{\text{a.s.}} 0,$$

*then a probability space can be constructed which supports a standard Wiener process $W$ and a process $Z'$ such that,*

$$\{Z(t) : t \ge 0\} \stackrel{\mathrm{d}}{=} \{Z'(t) : t \ge 0\} \quad and \quad \frac{|Z'(t) - W(t)|}{\sqrt{t[\log_2(t)]^{(1-\beta)}}} \xrightarrow{\text{a.s.}} 0.$$

*In the case where $\alpha = 1/2$, condition* (3.1) *may be replaced by the weaker condition*

$$(3.2) \qquad \left[ \frac{n}{[\log_2(n)]^{(1-\beta)}} \right]^{1/4} \left( \frac{R(n)}{n} - \alpha \right) \xrightarrow{\text{a.s.}} 0.$$

COROLLARY 3.1. *For any reading policy satisfying*

$$(3.3) \qquad \sqrt{n} \left( \frac{R(n)}{n} - \alpha \right) \xrightarrow{\text{a.s.}} 0, \qquad with \ \alpha \in (0,1),$$

*we have*

$$(3.4) \qquad \sqrt{n} \left[ \frac{M(n)}{\alpha(1-\alpha)n^2} - \mu \right] \xrightarrow{\mathrm{d}} \mathrm{N} \left( 0, \left[ \frac{(1-\alpha)\sigma_{\mathbf{R}}^2 + \alpha\sigma_{\mathbf{S}}^2}{\alpha(1-\alpha)} \right] \right).$$

*In the case where $\alpha = 1/2$, condition* (3.3) *may be replaced by*

$$(3.5) \qquad n^{1/4} \left( \frac{R(n)}{n} - \alpha \right) \xrightarrow{\text{a.s.}} 0.$$

REMARK 3.1. In Corollary 3.1 almost sure convergence in (3.3) and (3.5) may be replaced by convergence in probability. This follows from the proof of Theorem 3.1 and the martingale central limit theorem (see, e.g., Theorem 7.4 of [3]).

COROLLARY 3.2. *For any reading policy satisfying*

$$\sqrt{\frac{n}{\log_2 n}} \left( \frac{R(n)}{n} - \alpha \right) \xrightarrow{\text{a.s.}} 0, \tag{3.6}$$

*we have, with probability one,*

$$\limsup_{n \to \infty} \frac{M(n) - n^2 \alpha (1-\alpha)\mu}{n^{3/2} \sqrt{2\kappa \log_2 n}} = 1 = -\liminf_{n \to \infty} \frac{M(n) - n^2 \alpha (1-\alpha)\mu}{n^{3/2} \sqrt{2\kappa \log_2 n}}, \tag{3.7}$$

*where* $\kappa := \alpha(1-\alpha)((1-\alpha)\sigma_{\mathbf{R}}^2 + \alpha \sigma_{\mathbf{S}}^2)$. *When* $\alpha = 1/2$, *condition* (3.6) *may be replaced by*

$$\left[ \frac{n}{\log_2 n} \right]^{1/4} \left( \frac{R(n)}{n} - \alpha \right) \xrightarrow{\text{a.s.}} 0. \tag{3.8}$$

From the above corollaries we see that the central limit and iterated logarithm behavior of the number of matches are invariant among the class of policies whose selection ratios converge to $\alpha$ sufficiently fast. Included in this class are the alternating policies and (it will be shown) the greedy policies, both with $\alpha = 1/2$. Thus, these policies obey the same CLT and LIL, the optimality of the latter policy notwithstanding.

We observe that conditions (3.3) and (3.6) fail under Bernoulli sampling, where the source is determined by independent tosses of an $\alpha$-coin, $\alpha \neq 1/2$. Before discussing this case further, we give a simple example showing that these conditions can hold for a nondeterministic policy which imposes a restorative pressure to keep its selection ratio close to $\alpha$. Consider a reading policy with $R(1) = 1$ and which for $n > 1$ chooses source $\mathbf{R}$ with probability $\alpha_1 \in (0, \alpha)$ [resp., $\alpha_2 \in (\alpha, 1)$] when $R(n-1) \geq \alpha(n-1)$ [resp., $R(n-1) < \alpha(n-1)$]. For such a policy, we have $-U \overset{\text{st}}{\leq} R(n) - \alpha n \overset{\text{st}}{\leq} V$, where $U$ and $V$ are defined by

$$U := \min \left\{ k : \sum_{j=1}^{k} X_{2,j} \geq \alpha k + 1 \right\}$$

and

$$V := \min \left\{ k : \sum_{j=1}^{k} X_{1,j} \leq \alpha k - 1 \right\},$$

with $\{X_{i,j} : j \geq 1\}$ an i.i.d. $\mathrm{Ber}(\alpha_i)$ sequence for $i = 1, 2$. By Bernstein's inequality (e.g., see [15]), $U$ and $V$ have exponential tails from which conditions (3.3) and (3.6) follow.

We will see in the proof of Theorem 3.1 below that

$$\sqrt{n}\left[\frac{M(n)}{\alpha(1-\alpha)n^2} - \mu\right] \quad \text{and} \quad \sqrt{n}\left[\frac{(1-\alpha)\Gamma_{\mathbf{R}}[R(n)] + \alpha\Gamma_{\mathbf{S}}[S(n)]}{\alpha(1-\alpha)n} - 2\mu\right]$$

share the same weak limit quite generally and, in particular, under Bernoulli sampling. Moreover, note that under the Bernoulli sampling the expression $(1-\alpha)\Gamma_{\mathbf{R}}[R(n)] + \alpha\Gamma_{\mathbf{S}}[S(n)]$ is the $n$th partial sum of a sequence of independent variables, all with the distribution of $(1-\alpha)R(1)X_{\mathbf{R}}(1) + \alpha(1 - R(1))X_{\mathbf{S}}(1)$. Thus, by the ordinary CLT for i.i.d. sequences, we obtain the CLT (and, by a similar argument, the LIL) for $M(n)$ with the asymptotic variance given by

$$(3.9) \qquad \frac{\mu^2(1-2\alpha)^2}{\alpha(1-\alpha)} + \frac{(1-\alpha)\sigma_{\mathbf{R}}^2 + \alpha\sigma_{\mathbf{S}}^2}{\alpha(1-\alpha)}.$$

More generally, when $\{R(n)\}_{n \geq 1}$ is independent of the labels, the CLT holds under $\sqrt{n}(R(n)/n - \alpha) \xrightarrow{\mathrm{d}} \mathrm{N}(0, \sigma^2)$, with the asymptotic variance given by

$$\left[\frac{\mu(1-2\alpha)}{\alpha(1-\alpha)}\right]^2 \sigma^2 + \frac{(1-\alpha)\sigma_{\mathbf{R}}^2 + \alpha\sigma_{\mathbf{S}}^2}{\alpha(1-\alpha)}.$$

The argument uses Kolmogorov's maximal inequality to show that $Y_n$ in (3.12) is sufficiently close to the partial sum of the first $n\alpha$ $X_{\mathbf{R}}(\cdot)$'s and $n(1-\alpha)$ $X_{\mathbf{S}}(\cdot)$'s. The CLT follows from the independence of this partial sum and $R(n)$.

It is interesting to note the more stringent requirement in the above corollaries on the rate of convergence of the selection ratio when the limit is other than $1/2$. This is to account for the phenomenon that while the policy which uses Bernoulli sampling with $\alpha = 1/2$ obeys the same CLT as an alternating policy (and, as we shall see, a greedy policy), the policy which uses Bernoulli sampling with $\alpha \neq 1/2$ has a higher asymptotic variance than a policy for which $R(n) = \lceil n\alpha \rceil$ [cf. the asymptotic variance in (3.4) to the expression in (3.9) for the cases $\alpha = 1/2$ and $\alpha \neq 1/2$].

We now state the CLT and LIL for the selection ratio of a greedy policy, the latter result yielding both the CLT and LIL for the number of matches via an application of Corollary 3.1 and Corollary 3.2, respectively.

THEOREM 3.2. *For the greedy reading policy $C_{\mathbf{G}}$, we have*

$$(3.10) \qquad \frac{(R_{\mathbf{G}}(n) - n/2)}{\sqrt{n}} \xrightarrow{\mathrm{d}} \mathrm{N}(0, \sigma_{R_{\mathbf{G}}}^2) \qquad \textit{as } n \to \infty;\ \sigma_{R_{\mathbf{G}}}^2 := \left(\frac{\sigma_{\mathbf{R}}^2 + \sigma_{\mathbf{S}}^2}{8\mu^2}\right).$$

THEOREM 3.3. *For the greedy reading policy $C_{\mathbf{G}}$, we have*

$$(3.11) \quad \limsup_{n \to \infty} \frac{R_{\mathbf{G}}(n) - n/2}{\sqrt{2\sigma_{R_{\mathbf{G}}}^2 \, n \log_2 n}} = 1 = -\liminf_{n \to \infty} \frac{R_{\mathbf{G}}(n) - n/2}{\sqrt{2\sigma_{R_{\mathbf{G}}}^2 \, n \log_2 n}}, \qquad w.p. \ 1.$$

COROLLARY 3.3. *For both greedy and alternating policies, we have*

$$\sqrt{n}\left[\frac{M(n)}{(n/2)^2} - \mu\right] \xrightarrow{\mathrm{d}} \mathrm{N}(0, 2(\sigma_{\mathbf{R}}^2 + \sigma_{\mathbf{S}}^2)).$$

*Moreover, for these policies, we also have* (3.7) *with* $\alpha = 1/2$.

While both of the above theorems are of independent interest, the former is of interest also for the similarity of its derivation to that of the weak limit of a sequence of stopping times needed in the next section, and the latter for its application to the number of matches.

3.1. *Proofs of Theorems* 3.1–3.3.

PROOF OF THEOREM 3.1. First, we observe that

$$(3.12) \quad \begin{aligned} &\frac{M(n)}{n} - \mu_\alpha n - Y_n \\ &= \mu[1 - 2\alpha](R(n) - \alpha n) \\ &\quad + \alpha(1 - \alpha)n\left(\frac{\tilde{N}_{\mathbf{R}}(R(n))}{\alpha n} - \tilde{r}\right) \cdot \left(\frac{\tilde{N}_{\mathbf{S}}(S(n))}{(1 - \alpha)n} - \tilde{s}\right), \end{aligned}$$

where $\mu_\alpha := \alpha(1 - \alpha)\mu$ and for $n \geq 1$,

$$(3.13) \quad Y_n := (1 - \alpha)\Gamma_{\mathbf{R}}[R(n)] + \alpha\Gamma_{\mathbf{S}}[S(n)] - [\alpha n + (1 - 2\alpha)R(n)]\mu.$$

Second, we show that

$$(3.14) \quad \frac{(M(n))/n - \mu_\alpha n - Y_n}{\sqrt{n[\log_2(n)]^{(1-\beta)}}} \xrightarrow{\text{a.s.}} 0.$$

Let $a_n := \sqrt{n/[\log_2(n)]^{(1-\beta)}}$ for $n \geq 1$. Toward showing (3.14), we observe that when $\alpha = 1/2$, the term $a_n\mu[1 - 2\alpha](R(n)/n - \alpha)$ is 0 and otherwise it converges to 0 in the almost sure sense by (3.1). We note that this is the only reason for requiring the stringent condition (3.1). Now the proof of (3.14) is completed if we can show that (3.2) implies that

$$\left|a_n\left(\frac{\tilde{N}_{\mathbf{R}}(R(n))}{\alpha n} - \tilde{r}\right) \cdot \left(\frac{\tilde{N}_{\mathbf{S}}S(n)}{(1 - \alpha)n} - \tilde{s}\right)\right| \xrightarrow{\text{a.s.}} 0.$$

By the Cauchy–Schwarz inequality, we have

$$
\begin{aligned}
(3.15) \quad & \left| a_n \left( \frac{\tilde{N}_{\mathbf{R}} R(n)}{\alpha n} - \tilde{r} \right) \cdot \left( \frac{\tilde{N}_{\mathbf{S}}(S(n))}{(1-\alpha)n} - \tilde{s} \right) \right| \\
& \leq \sqrt{ \sum_{i=1}^{\infty} \left[ \sqrt{a_n} \left( \frac{N_{\mathbf{R}}[R(n), i]}{\alpha n} - r_i \right) \right]^2 } \sqrt{ \sum_{i=1}^{\infty} \left[ \sqrt{a_n} \left( \frac{N_{\mathbf{S}}[S(n), i]}{(1-\alpha)n} - s_i \right) \right]^2 }.
\end{aligned}
$$

By symmetry, it suffices to show that the first term on the right-hand side of (3.15) converges to zero in the almost sure sense:

$$
\sum_{i=1}^{\infty} \left[ \sqrt{a_n} \left( \frac{N_{\mathbf{R}}[R(n), i]}{\alpha n} - r_i \right) \right]^2
$$

$$
= \underbrace{ \frac{a_n \log_2 R(n)}{R(n)} \left( \frac{R(n)}{\alpha n} \right)^2 }_{\xrightarrow{\text{a.s.}} 0} \underbrace{ \frac{R(n)}{\log_2 R(n)} \sum_{i=1}^{\infty} \left( \frac{N_{\mathbf{R}}[R(n), i]}{R(n)} - r_i \right)^2 }_{O(1) \text{ a.s. by Lemma A.1}}
$$

$$
+ 2 \underbrace{ \sqrt{ \frac{a_n \log_2 R(n)}{R(n)} } \left( \frac{R(n)}{\alpha n} \right) }_{\xrightarrow{\text{a.s.}} 0} \underbrace{ \sqrt{a_n} \left( \frac{R(n)}{\alpha n} - 1 \right) }_{\xrightarrow{\text{a.s.}} 0 \text{ by (3.2)}}
$$

$$
\times \underbrace{ \sqrt{ \frac{R(n)}{\log_2 R(n)} } \sum_{i=1}^{\infty} \left( \frac{N_{\mathbf{R}}[R(n), i]}{R(n)} - r_i \right) r_i }_{O(1) \text{ a.s. by Lemma A.1}}
$$

$$
+ \underbrace{ a_n \left( \frac{R(n)}{\alpha n} - 1 \right)^2 \sum_{i=1}^{\infty} (r_i)^2 }_{\xrightarrow{\text{a.s.}} 0 \text{ by (3.2)}} \xrightarrow{\text{a.s.}} 0.
$$

Third, we show that $\{Y_n\}_{n \geq 1}$ is a $\{\mathcal{F}_n\}_{n \geq 1}$ martingale with bounded increments. Toward this, we note that, for $n \geq 2$,

$$
D_n := Y_n - Y_{n-1} = \begin{cases} (1-\alpha)(X_{\mathbf{R}}(R(n)) - \mu), & \text{if } C(n) = 1, \\ \alpha(X_{\mathbf{S}}(S(n)) - \mu), & \text{if } C(n) = 0, \end{cases}
$$

with $D_1 := Y_1$.

Now as $C(n)$ is $\mathcal{F}_{n-1}$ measurable and both $X_{\mathbf{R}}(R(n-1)+1)$ as well as $X_{\mathbf{S}}(S(n-1)+1)$ are independent of $\mathcal{F}_{n-1}$, we have $\mathbb{E}(D_n | \mathcal{F}_{n-1}) = 0$. Moreover, as $D_n$ is bounded and $Y_n$ is $\mathcal{F}_n$ measurable, we have the above description of $\{Y_n\}_{n \geq 1}$. Further, observe that $\mathbb{E}(D_n^2 | \mathcal{F}_{n-1}) = (1-\alpha)^2 \sigma_R^2 C(n) + \alpha^2 \sigma_S^2 (1 - C(n))$, which implies that

$$
\frac{1}{n} \sum_{k=1}^{n} \mathbb{E}(D_k^2 | \mathcal{F}_{k-1}) = (1-\alpha)^2 \sigma_R^2 \frac{R(n)}{n} + \alpha^2 \sigma_S^2 \frac{S(n)}{n}
$$

$$\xrightarrow{\text{a.s.}} \alpha(1-\alpha)((1-\alpha)\sigma_{\mathbf{R}}^2 + \alpha\sigma_{\mathbf{S}}^2).$$

Fourth, we observe that the above description of $\{Y_n\}_{n\geq 1}$ along with Theorem 3.2 of [8] implies that a probability space can be constructed (a suitably augmented version of the one in [8]) which supports a Wiener process $W$ and a sequence $\{(R'(n), L'_{\mathbf{R}}(n), L'_{\mathbf{S}}(n))\}_{n\geq 1}$ satisfying the following:

(i) $\{(R'(n), L'_{\mathbf{R}}(n), L'_{\mathbf{S}}(n))\}_{n\geq 1} \overset{\mathrm{d}}{=} \{(R(n), L_{\mathbf{R}}(n), L_{\mathbf{S}}(n))\}_{n\geq 1}$,

(ii) If $Y'_n$ is the same function of $\{(R'(n), L'_{\mathbf{R}}(n), L'_{\mathbf{S}}(n))\}_{n\geq 1}$ as $Y_n$ is of $\{(R(n), L_{\mathbf{R}}(n), L_{\mathbf{S}}(n))\}_{n\geq 1}$ and

$$Y'(t) := \begin{cases} 0, & \text{for } 0 \leq t < V'_1, \\ Y'(n), & \text{for } V'_n \leq t < V'_{n+1} \text{ and } n \geq 1, \end{cases}$$

where $V'_n := (1-\alpha)^2\sigma_R^2 R'(n) + \alpha^2\sigma_S^2 S'(n)$ for $n \geq 1$, then

(3.16)
$$\frac{|Y'(t) - W(t)|}{\sqrt{t[\log_2(t)]^{(1-\beta)}}} \xrightarrow{\text{a.s.}} 0.$$

To complete the proof, let $\{M'(n)\}_{n\geq 1}$ be the sequence defined exactly as $\{M(n)\}_{n\geq 1}$ but using the sequence $\{(R'(n), L'_{\mathbf{R}}(n), L'_{\mathbf{S}}(n))\}_{n\geq 1}$ instead of $\{(R(n), L_{\mathbf{R}}(n), L_{\mathbf{S}}(n))\}_{n\geq 1}$. Also, let $Z'$ be the process defined like $Z$ but using $\{(M'(n), V'(n))\}_{n\geq 1}$ instead of $\{(M(n), V(n))\}_{n\geq 1}$. Then (3.14) and (3.16) together imply

$$\frac{|Z'(t) - W(t)|}{\sqrt{t[\log_2(t)]^{(1-\beta)}}} \xrightarrow{\text{a.s.}} 0.$$

Hence, the proof.  □

To prove Theorems 3.2 and 3.3, we utilize the following lemma. This lemma provides a tight connection with partial sums of i.i.d. variables that is key in the proofs of the two theorems.

LEMMA 3.1.   *In the case of the greedy algorithm we have, for* $0 \leq x \leq n$,

$$\Gamma_{\mathbf{R}}[\lceil x \rceil] < \Gamma_{\mathbf{S}}[n - \lceil x \rceil] - \gamma \quad \Longrightarrow \quad R_{\mathbf{G}}(n) > x$$
$$\Longrightarrow \quad \Gamma_{\mathbf{R}}[\lfloor x \rfloor] \leq \Gamma_{\mathbf{S}}[n - \lfloor x \rfloor] + \gamma.$$

PROOF.   Since $\Gamma_{\mathbf{S}}[\cdot]$ and $\Gamma_{\mathbf{R}}[\cdot]$ are nondecreasing, for $0 \leq x \leq n$, we have, by (2.5),

$$R_{\mathbf{G}}(n) > x$$
$$\Longrightarrow \quad \Gamma_{\mathbf{S}}[\lfloor x \rfloor] \leq \Gamma_{\mathbf{R}}[R_{\mathbf{G}}(n)] \leq \Gamma_{\mathbf{S}}[S_{\mathbf{G}}(n)] + \gamma \leq \Gamma_{\mathbf{S}}[n - \lfloor x \rfloor] + \gamma.$$

The other half follows by observing that, for $0 \le x \le n$,

$$R_{\mathbf{G}}(n) \le x$$
$$\implies \Gamma_{\mathbf{R}}[\lceil x \rceil] \ge \Gamma_{\mathbf{R}}[R_{\mathbf{G}}(n)] \ge \Gamma_{\mathbf{S}}[S_{\mathbf{G}}(n)] - \gamma \ge \Gamma_{\mathbf{S}}[n - \lceil x \rceil] - \gamma. \quad \square$$

PROOF OF THEOREM 3.2. First, we argue that

$$
\begin{aligned}
(3.17) \quad Z(n) &:= \frac{\Gamma_{\mathbf{R}}[k_n] - \Gamma_{\mathbf{S}}[n - k_n]}{\sqrt{n/2}} \\
&\xrightarrow{\mathrm{d}} \mathrm{N}(2^{3/2}\mu x, \sigma_{\mathbf{R}}^2 + \sigma_{\mathbf{S}}^2), \qquad \text{as } n \to \infty
\end{aligned}
$$

for any sequence $\{k_n\}_{n \ge 1}$ satisfying

$$\lim_{n \to \infty} \left( \frac{k_n - n/2}{\sqrt{n}} \right) = x.$$

To this end, note that

$$
\begin{aligned}
(3.18) \quad Z(n) = a_n &\underbrace{\left[ \frac{\sum_{j=1}^{k_n} X_{\mathbf{R}}(j) - k_n \mu}{\sqrt{k_n}} \right]}_{\xrightarrow{\mathrm{d}} \mathrm{N}(0, \sigma_{\mathbf{R}}^2)} \\
&- b_n \underbrace{\left[ \frac{\sum_{j=1}^{n-k_n} X_{\mathbf{S}}(j) - (n - k_n)\mu}{\sqrt{n - k_n}} \right]}_{\xrightarrow{\mathrm{d}} \mathrm{N}(0, \sigma_{\mathbf{S}}^2)} + c_n 2^{3/2} \mu x,
\end{aligned}
$$

where the three sequence $\{a_n\}_{n \ge 1}$, $\{b_n\}_{n \ge 1}$ and $\{c_n\}_{n \ge 1}$ all converge to 1. The above observed weak limits are due to the ordinary central limit theorem. By independence of the first two terms in (3.18) and Slutsky's theorem, we have (3.17). Now defining $Z_*$ and $Z^*$ as $Z$ but with the sequence $\{k_n\}_{n \ge 1}$ taken as $\{\lceil n/2 + x\sqrt{n} \rceil\}_{n \ge 1}$ and $\{\lfloor n/2 + x\sqrt{n} \rfloor\}_{n \ge 1}$, respectively, we have, as $n \to \infty$,

$$
\begin{aligned}
(3.19) \quad Z_*(n) &\xrightarrow{\mathrm{d}} \mathrm{N}(2^{3/2}\mu x, \sigma_{\mathbf{R}}^2 + \sigma_{\mathbf{S}}^2) \quad \text{and} \\
Z^*(n) &\xrightarrow{\mathrm{d}} \mathrm{N}(2^{3/2}\mu x, \sigma_{\mathbf{R}}^2 + \sigma_{\mathbf{S}}^2).
\end{aligned}
$$

By Lemma 3.1, we have for large $n$

$$
\begin{aligned}
(3.20) \quad \Pr\left( Z_*(n) < \frac{-\gamma}{\sqrt{n/2}} \right) &\le \Pr\left( \frac{R_{\mathbf{G}} n - n/2}{\sqrt{n}} > x \right) \\
&\le \Pr\left( Z^*(n) < \frac{\gamma}{\sqrt{n/2}} \right),
\end{aligned}
$$

which combined with (3.19) completes the proof. $\square$

PROOF OF THEOREM 3.3.   Let $\{k_n\}_{n\geq 1}$ be a sequence of nonnegative integers and $\{a_n\}_{n\geq 1}$ a sequence of reals such that, for $n \geq 1$,

$$\frac{k_n - n/2}{\sqrt{2\sigma_{R_\mathbf{G}}^2 n \log_2 n}} \xrightarrow{n \to \infty} C > 0 \quad \text{and}$$

(3.21)

$$a_n = \frac{1}{\sqrt{2(\sigma_\mathbf{R}^2 + \sigma_\mathbf{S}^2)(n - k_n)\log_2(n - k_n)}}.$$

For such a sequence $\{k_n\}_{n\geq 1}$,

$$a_n(\Gamma_\mathbf{R}[k_n] - \Gamma_\mathbf{S}[n - k_n]) = \underbrace{a_n(\Gamma_\mathbf{R}[n - k_n] - \Gamma_\mathbf{S}[n - k_n])}_{\liminf_{n\to\infty} = -1}$$

$$+ \underbrace{a_n(\Gamma_\mathbf{R}[k_n] - \Gamma_\mathbf{R}[n - k_n] - (2k_n - n)\mu)}_{\xrightarrow{\text{a.s.}} 0} + \underbrace{b_n}_{\xrightarrow{n \to \infty} C},$$

where the first limit infimum is due to the standard law of iterated logarithm, the second limit is due to Theorem 5.1 of [6] on lag sums and the third limit is a consequence of (3.21). Hence, for a sequence $\{k_n\}_{n\geq 1}$ satisfying (3.21), we have

$$(3.22) \qquad \liminf_{n \to \infty} a_n(\Gamma_\mathbf{R}[k_n] - \Gamma_\mathbf{S}[n - k_n]) = C - 1.$$

Using (3.22) and Lemma 3.1, with $k_n = \lceil n/2 + (1 - \varepsilon)\sqrt{2\sigma_{R_\mathbf{G}}^2 n \log_2 n} \rceil$, we have

$$(3.23) \quad \frac{R_\mathbf{G}(n) - n/2}{\sqrt{2\sigma_{R_\mathbf{G}}^2 n \log_2 n}} > 1 - \varepsilon \qquad \text{infinitely often (i.o.) a.s. } \forall \varepsilon > 0.$$

Now similarly, working instead with $k_n = \lfloor n/2 + (1 + \varepsilon)\sqrt{2\sigma_{R_\mathbf{G}}^2 n \log_2 n} \rfloor$, we have

$$(3.24) \quad \frac{R_\mathbf{G}(n) - n/2}{\sqrt{2\sigma_{R_\mathbf{G}}^2 n \log_2 n}} > 1 + \varepsilon \qquad \text{only finitely often a.s. } \forall \varepsilon > 0.$$

Statements (3.23) and (3.24) are equivalent to the first statement in (3.11). A similar argument leads to the other. Hence, the proof. $\square$

**4. Comparison of greedy and the alternating.** The results of the last section, which say that the weak limit and the law of the iterated logarithm for $M_\mathbf{G}(\cdot)$ and $M_\mathbf{A}(\cdot)$ coincide, motivate asymptotic analysis of their difference—the goal of this section.

In our problem the reward is unbounded and grows linearly. Hence, an analogue to the average reward criterion in the bounded reward case would involve $\mathbb{E}(M(n))/n^2$. We note that by Corollary 3.3 and the dominated convergence theorem we have the equality of $\lim_{n\to\infty} n^{-2}\mathbb{E}(M_{\mathbf{G}}(n))$ and $\lim_{n\to\infty} n^{-2}\mathbb{E}(M_{\mathbf{A}}(n))$. This leads us to consider sensitive discount optimality criteria to distinguish the performance of the alternating from that of the greedy. Let us denote by $\nu_\lambda^{\mathrm{A}}$ (resp., $\nu_\lambda^{\mathrm{G}}$), for $0 \le \lambda < 1$, the expected total $\lambda$-discounted incremental matches under the alternating policy (resp., greedy policy). In the case of the alternating policy an easy calculation yields $\nu_\lambda^{\mathrm{A}} = (1-\lambda)^{-2}\lambda(1+\lambda)^{-1}\mu$. A study of the $-1$-discount optimality of the alternating leads us to $\liminf_{\lambda\uparrow 1}(1-\lambda)(\nu_\lambda^{\mathrm{G}} - \nu_\lambda^{\mathrm{A}})$. It follows by a Tauberian theorem of Hardy and Littlewood (see, e.g., Theorem 7.4 of [9]) that

$$\lim_{\lambda\uparrow 1}(1-\lambda)(\nu_\lambda^{\mathrm{G}} - \nu_\lambda^{\mathrm{A}}) = \lim_{n\to\infty}\left(\frac{\mathbb{E}(M_{\mathbf{G}}(n) - M_{\mathbf{A}}(n))}{n}\right)$$

when either limit exists. The first theorem shows that the limit on the right exists and is positive, hence, showing that the alternating policy is not $-1$-discount optimal.

THEOREM 4.1. *For two chosen policies, one greedy and the other alternating, we have*

$$\lim_{n\to\infty}\left(\frac{\mathbb{E}(M_{\mathbf{G}}(n) - M_{\mathbf{A}}(n))}{n}\right) = \frac{\sigma_{\mathbf{R}}^2 + \sigma_{\mathbf{S}}^2}{8\mu}.$$

REMARK 4.1. It is easily seen that $\mathbb{E}(\max(\Gamma_{\mathbf{R}}[R_{\mathbf{G}}(n)], \Gamma_{\mathbf{S}}[S_{\mathbf{G}}(n)]))$ represents the expected incremental gain of matches from the $(n+1)$st pick by a greedy policy. The proof of Theorem 3.1 then gives us

$$\mathbb{E}(\max(\Gamma_{\mathbf{R}}[R_{\mathbf{G}}(n)], \Gamma_{\mathbf{S}}[S_{\mathbf{G}}(n)])) = \left(\frac{1}{2}\right)\mathbb{E}(|\Gamma_{\mathbf{R}}[R_{\mathbf{G}}(n)] - \Gamma_{\mathbf{S}}[S_{\mathbf{G}}(n)]|)$$

$$+ \left(\frac{1}{2}\right)\mathbb{E}(\Gamma_{\mathbf{R}}[R_{\mathbf{G}}(n)] + \Gamma_{\mathbf{S}}[S_{\mathbf{G}}(n)])$$

$$= \left(\frac{1}{2}\right)\mathbb{E}(|\Gamma_{\mathbf{R}}[R_{\mathbf{G}}(n)] - \Gamma_{\mathbf{S}}[S_{\mathbf{G}}(n)]|) + \frac{n\mu}{2}.$$

Interestingly, the greedy criterion implies that the process $\{\Gamma_{\mathbf{R}}[R_{\mathbf{G}}(n)] - \Gamma_{\mathbf{S}}[S_{\mathbf{G}}(n)]\}_{n\ge 1}$ is a bounded Markov chain on a subset of $[-\gamma, \gamma]$, leaving aside the versions of greedy which introduce unnecessary path dependence on the epochs where the greedy criterion is ambivalent. This immediately leads to the relation

$$(4.1) \quad \mathbb{E}(M_{\mathbf{G}}(n)) = \left(\frac{1}{2}\right)\sum_{k=1}^{n-1}\mathbb{E}(|\Gamma_{\mathbf{R}}[R_{\mathbf{G}}(n)] - \Gamma_{\mathbf{S}}[S_{\mathbf{G}}(n)]|) + \frac{n(n-1)\mu}{4}.$$

REMARK 4.2. In the case of ergodicity of $\{\Gamma_{\mathbf{R}}[R_{\mathbf{G}}(n)] - \Gamma_{\mathbf{S}}[S_{\mathbf{G}}(n)]\}_{n \geq 1}$, the above theorem yields

$$\mathbb{E}(\Delta M_{\mathbf{G}}(n) - \Delta M_{\mathbf{A}}(n)) \to \begin{cases} \dfrac{\sigma_{\mathbf{R}}^2 + \sigma_{\mathbf{S}}^2}{8\mu} - \dfrac{\mu}{4}, & \text{along odd } n\text{'s}, \\ \dfrac{\sigma_{\mathbf{R}}^2 + \sigma_{\mathbf{S}}^2}{8\mu} + \dfrac{\mu}{4}, & \text{along even } n\text{'s}, \end{cases}$$

where $\Delta$ is the difference operator. It is easily checked that there are examples where $\mathbb{E}(\Delta M_{\mathbf{G}}(n) - \Delta M_{\mathbf{A}}(n))$ for all large $n$ is positive and where it oscillates in sign.

REMARK 4.3. In the case of geometric ergodicity of $\{\Gamma_{\mathbf{R}}[R_{\mathbf{G}}(n)] - \Gamma_{\mathbf{S}}[S_{\mathbf{G}}(n)]\}_{n \geq 1}$, the above theorem together with (4.1) leads to an expansion of the form $\{\Gamma_{\mathbf{R}}[R_{\mathbf{G}}(n)] - \Gamma_{\mathbf{S}}[S_{\mathbf{G}}(n)]\}_{n \geq 1}$,

$$\mathbb{E}(M_{\mathbf{G}}(n)) = \frac{n^2 \mu}{4} + n \left[ \frac{\sigma_{\mathbf{R}}^2 + \sigma_{\mathbf{S}}^2}{8\mu} \right] + \text{constant} + \varepsilon_n,$$

where $\varepsilon_n$ tends to zero exponentially fast. This relation for the illustrative example is given in (1.1), and Figure 1 graphically describes the Markov chain $\{\Gamma_{\mathbf{R}}[R_{\mathbf{G}}(n)] - \Gamma_{\mathbf{S}}[S_{\mathbf{G}}(n)]\}_{n \geq 1}$.

REMARK 4.4. The above, in particular, implies that, for large $n$,

$$\mathbb{E}(M_{\mathbf{A}}(n)) \leq \mathbb{E}(M_{\mathbf{G}}(n)) \leq \mathbb{E}(M_{\mathbf{A}}(n + k)), \qquad \text{where } k := \left\lceil \frac{\sigma_{\mathbf{R}}^2 + \sigma_{\mathbf{S}}^2}{4\mu^2} \right\rceil.$$

An easy but informative upper bound for the $k$ of the above equation is given by $\lceil [(1 - \mu)/2\mu] \rceil$. It can be easily shown that, in general, $k$ cannot be bounded away from infinity.

The first theorem, while able to distinguish between the greedy and the alternating policies, also shows that the difference is rather slim. It would be interesting to know whether these two policies, if implemented on the same sampled sequence of labels, will yield (in some sense) more matches under the greedy than the alternating? Such optimality is referred to in the literature as sample path optimality. Results on sample path optimality for the case of bounded rewards and finite/countable state and action spaces under the assumption of uniform ergodicity of the state process can be found in [7] and the references therein. The second theorem shows that the weak limit of $M_{\mathbf{G}}(n) - M_{\mathbf{A}}(n)$ under this coupling is a scale mixture of normals centered at zero. The final theorem of this section throws some light on its sample path behavior under the same coupling. A more precise study of the sample path behavior of $M_{\mathbf{G}}(n) - M_{\mathbf{A}}(n)$ is beyond the scope of this paper.

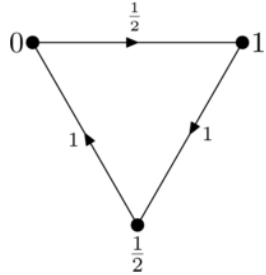

FIG. 1. *Embedded Markov chain of the illustrative example.*

THEOREM 4.2. *For two chosen policies, one greedy and the other alternating, we have*

$$\left(\frac{M_{\mathbf{G}}(n) - M_{\mathbf{A}}(n)}{n^{5/4}}\right) \xrightarrow{\mathrm{d}} \mathrm{F} \qquad as \ n \to \infty,$$

*where* $\mathrm{F}$ *is a scale mixture of normals centered at zero given by*

$$\mathrm{F} = \int \mathrm{N}(0, \sigma^2) \, dG(\sigma^2), \qquad where \ \mathrm{G} := \left| \mathrm{N}\left(0, \frac{(\sigma_{\mathbf{R}}^2 + \sigma_{\mathbf{S}}^2)^3}{128\mu^2}\right) \right|.$$

THEOREM 4.3. *For two chosen policies, one greedy and the other alternating, we have*

(4.2)
$$\liminf_{n \to \infty} \left(\frac{M_{\mathbf{G}}(n) - M_{\mathbf{A}}(n)}{n^{5/4}(\log_2 n)^{1/4}}\right) = -\infty \qquad a.s.,$$

$$\limsup_{n \to \infty} \left(\frac{M_{\mathbf{G}}(n) - M_{\mathbf{A}}(n)}{n^{5/4}(\log_2 n)^{1/4}}\right) = \infty \qquad a.s.$$

*and*

(4.3)
$$\frac{M_{\mathbf{G}}(n) - M_{\mathbf{A}}(n)}{n^{5/4}(\log_2 n)^{1/4}} = O((\log n)^{1/2}) \qquad a.s.$$

In the following we will need the filtration $\{\mathcal{G}_n\}_{n \geq 0}$ defined as

$$\mathcal{G}_n = \mathcal{G}_0 \vee \sigma\langle L_{\mathbf{R}}(1), \dots, L_{\mathbf{R}}(R_{\mathbf{G}}(n)); L_{\mathbf{S}}(1), \dots, L_{\mathbf{S}}(S_{\mathbf{G}}(n))\rangle, \qquad n \geq 1,$$

with $\mathcal{G}_0$ containing all the information needed for randomization by not only $C_{\mathbf{G}}$ but also $C_{\mathbf{A}}$. The argument for the above results depends on the sequence of random times $\{T_n\}_{n \geq 1}$, where $T_n$ is essentially the epoch at which the greedy decides to pick the first record (from $\mathbf{R}$ or $\mathbf{S}$) which would not be seen following the alternating policy by the $n$th epoch. Formally, they are defined as

$$T_n := \min\left(n, \inf\left\{k \geq 1 \,\Big|\, S_{\mathbf{G}}(k+1) = \left\lceil \frac{n}{2} \right\rceil + 1 \text{ or } R_{\mathbf{G}}(k+1) = \left\lceil \frac{n}{2} \right\rceil + 1\right\}\right)$$
(4.4)

for $n \geq 1$. It is easily checked that $\{T_n\}_{n \geq 1}$ is a sequence of $\{\mathcal{G}_n\}_{n \geq 0}$ stopping times. Also, for convenience, we define the sequence of events $\{A_n\}_{n \geq 1}$ by $A_n := \{R_{\mathbf{G}}(T_n) = \lceil n/2 \rceil\}$ for $n \geq 1$.

4.1. *Study of stopping times* $\{T_n\}_{n \geq 1}$. We note that, similarly to Lemma 3.1, it can be shown that, for positive $x$,

$$(n - T_n) > x \quad \Longrightarrow \quad \Gamma_{\mathbf{R}}[\lceil n/2 \rceil] \leq \Gamma_{\mathbf{S}}[\lfloor n/2 \rfloor - \lfloor x \rfloor] + \gamma \quad \text{or}$$
(4.5)
$$\Gamma_{\mathbf{S}}[\lceil n/2 \rceil] \leq \Gamma_{\mathbf{R}}[\lfloor n/2 \rfloor - \lfloor x \rfloor] + \gamma$$

and

$$(n - T_n) > x \quad \Longleftarrow \quad \Gamma_{\mathbf{R}}[\lceil n/2 \rceil] < \Gamma_{\mathbf{S}}[\lfloor n/2 \rfloor - \lceil x \rceil] \quad \text{or}$$
(4.6)
$$\Gamma_{\mathbf{S}}[\lceil n/2 \rceil] < \Gamma_{\mathbf{R}}[\lfloor n/2 \rfloor - \lceil x \rceil].$$

This leads to the first lemma which describes both the weak limit and sample path behavior of $T_n$. The second lemma is the weak law of large numbers for the post $T_n$ (and pre-$n$) selection ratio. The third lemma derives exponential probability inequalities for both $T_n$ and $R_{\mathbf{G}}(n)$ which are useful in establishing the required uniform integrability results and the uniform central limit theorem of Lemma 4.6.

LEMMA 4.1. *For the above defined stopping times* $\{T_n\}_{n \geq 1}$, *the following hold:*

    (i)

$$\frac{n - T_n}{2\sqrt{n}} \xrightarrow{\text{d}} |\mathrm{N}(0, \sigma_{R_{\mathbf{G}}}^2)| \qquad \text{as } n \to \infty,$$
(4.7)

*where* $\sigma_{R_{\mathbf{G}}}^2$, *the asymptotic variance of* $R_{\mathbf{G}}(n)$, *is defined in* (3.10).

    (ii)

$$\liminf (n - T_n) = 0 \qquad a.s. \quad and \quad \limsup \frac{n - T_n}{\sqrt{8\sigma_{R_{\mathbf{G}}}^2 \, n \log_2 n}} = 1 \qquad a.s.$$
(4.8)

PROOF. A proof of (4.7) and the second part of (4.8) follows along similar lines as Theorem 3.2 and Theorem 3.3, respectively. The key difference being that (4.5) and (4.6) are used instead of Lemma 3.1. The first part of (4.8) follows as Theorem 3.3 implies that the event $\{R_{\mathbf{G}}(n) = S_{\mathbf{G}}(n)\}$ for infinitey many $n\}$ occurs with probability one. The details are skipped to avoid repetition of similar arguments. □

LEMMA 4.2. *For the above defined stopping times $\{T_n\}_{n\geq 1}$ corresponding to the greedy reading policy $C_{\mathbf{G}}$, we have*

$$(4.9) \qquad \frac{R_{\mathbf{G}}(n) - R_{\mathbf{G}}(T_n)}{n - T_n} \xrightarrow{\ \mathrm{P}\ } \frac{1}{2} \qquad as\ n \to \infty.$$

PROOF. First, we show that

$$(4.10) \qquad \frac{R_{\mathbf{G}}(n) - R_{\mathbf{G}}(T_n)}{\log(n)} \xrightarrow{\ \mathrm{P}\ } \infty \qquad as\ n \to \infty.$$

By a double application of (2.5), we have

$$(4.11) \quad |(\Gamma_{\mathbf{R}}[R_{\mathbf{G}}(n)] - \Gamma_{\mathbf{R}}[R_{\mathbf{G}}(T_n)]) - (\Gamma_{\mathbf{S}}[S_{\mathbf{G}}(n)] - \Gamma_{\mathbf{S}}[S_{\mathbf{G}}(T_n)])| \leq 2\gamma,$$

which implies that, for any positive $K$,

$$(4.12) \quad R_{\mathbf{G}}(n) - R_{\mathbf{G}}(T_n) < K\log(n) \implies$$
$$\sum_{i=1}^{n-T_n-K\log(n)} X_{\mathbf{S}}(i + S_{\mathbf{G}}(T_n)) - 2\gamma < \sum_{i=1}^{K\log(n)} X_{\mathbf{R}}(i + R_{\mathbf{G}}(T_n)).$$

The second expression can be rewritten as

$$\left( \sum_{i=1}^{n-T_n-K\log(n)} \frac{X_{\mathbf{S}}(i + S_{\mathbf{G}}(T_n)) - \mu}{\sqrt{n - T_n - K\log(n)}} \right) - \left( \sum_{i=1}^{K\log(n)} \frac{X_{\mathbf{R}}(i + R_{\mathbf{G}}(T_n)) - \mu}{\sqrt{K\log(n)}} \right)$$
$$< 2\gamma - \mu(\sqrt{n - T_n - K\log(n)} - \sqrt{K\log(n)}).$$

As $n^{-1/4}(n - T_n) \xrightarrow{\ \mathrm{P}\ } \infty$ (Lemma 4.1), we have the independent terms on the left converging to normal distributions and the term on the right converging to negative infinity. Hence, the probability of the above event converges to zero. Now using (4.12), we obtain (4.10). Combining (4.10) with Theorem 3.3, we have

$$(4.13) \quad R_{\mathbf{G}}(n) \xrightarrow{\ \mathrm{P}\ } \infty \quad \text{and} \quad \frac{R_{\mathbf{G}}(n) - R_{\mathbf{G}}(T_n)}{\log(R_{\mathbf{G}}(n))} \xrightarrow{\ \mathrm{P}\ } \infty \qquad as\ n \to \infty.$$

Statement (4.13) together with Theorem 5.1 of [6] on lag sums gives us

$$\left( \frac{\Gamma_{\mathbf{R}}[R_{\mathbf{G}}(n)] - \Gamma_{\mathbf{R}}[R_{\mathbf{G}}(T_n)]}{R_{\mathbf{G}}(n) - R_{\mathbf{G}}(T_n)} \right) \xrightarrow{\ \mathrm{P}\ } \mu$$

and

$$\left( \frac{\Gamma_{\mathbf{S}}[S_{\mathbf{G}}(n)] - \Gamma_{\mathbf{S}}[S_{\mathbf{G}}(T_n)]}{S_{\mathbf{G}}(n) - S_{\mathbf{G}}(T_n)} \right) \xrightarrow{\ \mathrm{P}\ } \mu,$$

where the second part follows by symmetry. This with (4.11) and (4.10) gives

$$\frac{R_{\mathbf{G}}(n) - R_{\mathbf{G}}(T_n)}{S_{\mathbf{G}}(n) - S_{\mathbf{G}}(T_n)} \xrightarrow{\mathrm{P}} 1,$$

which is equivalent to (4.9). $\square$

REMARK 4.5. The above, in particular, implies that

$$\frac{R_{\mathbf{G}}(n) - R_{\mathbf{G}}(T_n)}{\sqrt{n}} = \left(\frac{R_{\mathbf{G}}(n) - R_{\mathbf{G}}(T_n)}{n - T_n}\right)\left(\frac{n - T_n}{\sqrt{n}}\right)$$

$$\xrightarrow{\mathrm{d}} |\,\mathrm{N}(0, \sigma^2_{R_{\mathbf{G}}})| \qquad \text{as } n \to \infty,$$

where the convergence in probability of the first term follows from (4.9) and the weak convergence of the second term to the folded normal follows from (4.7). In fact, the stronger result

$$\frac{R_{\mathbf{G}}(n) - R_{\mathbf{G}}(T_n) - (n - T_n)/2}{\sqrt{n - T_n}} \xrightarrow{\mathrm{d}} \mathrm{N}(0, \sigma^2_{R_{\mathbf{G}}}) \qquad \text{as } n \to \infty$$

can be shown using an argument similar to that used to prove Theorem 3.2.

LEMMA 4.3. *For a greedy policy, we have the following:*

(i) *For $t \geq 1$ and $n \geq (2 + \frac{\gamma}{\mu})^2$,*

$$(4.14) \qquad \Pr\left(\left|\frac{R_{\mathbf{G}}(n) - n/2}{\sqrt{n}}\right| > t\right) \leq 2 \exp\left\{-\left(\frac{\mu}{2\gamma}\right)^2 t\right\}.$$

(ii) *The sequence*

$$(4.15) \qquad \left\{\left(\frac{n - T_n}{\sqrt{n}}\right)^2\right\}_{n \geq 1}$$

*is uniformly integrable.*

PROOF. Toward proving (4.14), let $k = \lfloor n/2 + t\sqrt{n} \rfloor$. By Lemma 3.1, we have

$$\Pr(R_{\mathbf{G}} n > n/2 + t\sqrt{n}) \leq \Pr(\Gamma_{\mathbf{R}}[k] \leq \Gamma_{\mathbf{S}}[n - k] + \gamma).$$

Observe that

$$\Gamma_{\mathbf{S}}[n - k] - \Gamma_{\mathbf{R}}[k] + (2k - n)\mu = [\Gamma_{\mathbf{S}}[n - k] - \Gamma_{\mathbf{R}}[n - k]]$$

$$- [\Gamma_{\mathbf{R}}[k] - \Gamma_{\mathbf{R}}[n - k] - (2k - n)\mu],$$

which implies that the term on the left-hand side is a sum of $k$ independent zero mean random variables taking values in $[-\gamma, \gamma]$. This along with Hoeffding's inequality (see [15]) implies

$$\Pr \Gamma_{\mathbf{R}}[k] \leq \Gamma_{\mathbf{S}}[n-k] + \gamma \leq \exp\left\{\frac{-[(2k-n)\mu - \gamma]^2}{2k\gamma^2}\right\}.$$

Working with the upper bound above and using the inherent symmetry, we get the simple upper bound in (4.14).

For the sequence in (4.15), we get an inequality similar to (4.14) by imitating the above argument—the only change being that (4.5) and (4.6) are used instead of Lemma 3.1. Since this bound is integrable and free of $n$, we have the uniform integrability of the sequence in (4.15).    □

4.2. *Heuristics for the theorems.* We find it convenient to divide the records sampled by the $n$th epoch by either the alternating or the greedy into three sets. The first set consists of the first $T_n$ records sampled by the greedy. The second consists of the last $n - T_n$ records sampled by the greedy. The third consists of the records sampled by the alternating and not contained in the first set. Observe that all of the records in the first set (except possibly one) are sampled by the alternating by the $n$th epoch. Also note that all of the records in the third set belong to a single source and its cardinality is within one of the cardinality of the second set. The upshot of this is that $M_{\mathbf{G}}(n) - M_{\mathbf{A}}(n)$ is essentially the number of matches between records of the first and the second set and the records of the second set with themselves minus the number of matches between records of the first and the third set.

First, we argue that in the expected difference $\mathbb{E}(M_{\mathbf{G}}(n) - M_{\mathbf{A}}(n))$ the significant term comes from the matches generated among the records of the second set, which by Lemmas 4.1 and 4.3 will be of order $n$. This is so as the expected number of matches between records of the first and the second set minus the expected number of matches between the first and third set is at the most of order $\sqrt{n}$—follows by observing that $|\Gamma_{\mathbf{R}}[R_{\mathbf{G}}(T_n)] - \Gamma_{\mathbf{S}}[S_{\mathbf{G}}(T_n)]| \leq \gamma$ and $n - T_n = O_p(\sqrt{n})$ (by Lemma 4.1). Second, by Lemma 4.2, roughly $(n - T_n)/2$ of the records in the second set will be from each source and, hence, using the law of large numbers, the expected number of matches generated by these among themselves will be approximately $\mathbb{E}((n-T_n)^2)\mu/4$ or by Lemma 4.1, approximately $n\mu\sigma_{R_{\mathbf{G}}}^2$. This completes the heuristic for Theorem 4.1. Lemma 4.4 formalizes the latter part of the argument and the proof of Theorem 4.1 does the rest.

In contrast to the above, in the study of the sample path behavior and weak limit of $M_{\mathbf{G}}(n) - M_{\mathbf{A}}(n)$ we find that the insignificant terms of the above become significant and vice-versa. First, the number of matches generated by the records of the second set with themselves is comparable to $M_{\mathbf{G}}(\sqrt{n - T_n})$ which is $O_p(n)$, using Corollary 3.1 and Lemma 4.1. On the other hand, the number of matches between records of the first and second set minus the number of matches between the first and third set is $O_p(n^{5/4})$. Toward an argument, suppose without loss of generality that the greedy picks more $\mathbf{R}$ records than the alternating. Now the above difference is easily checked to be the difference in the numbers of matches generated by the *excess* $\mathbf{R}$ records sampled by the greedy and the *excess* $\mathbf{S}$ records sampled by the alternating, both with records from the first set. And divided by $n$, by the law of large numbers, the distribution of labels on the records from both sources in the first set approaches their respective vectors ($\tilde{r}$ or $\tilde{s}$). This then makes the difference resemble $nG_n$, where the sequence of random variables $\{G_n\}_{n \geq 1}$ is defined as

$$2G_n := \begin{cases} (\Gamma_{\mathbf{R}}[R_{\mathbf{G}}(n)] - \Gamma_{\mathbf{R}}[\lceil n/2 \rceil]) - (\Gamma_{\mathbf{S}}[\lfloor n/2 \rfloor] - \Gamma_{\mathbf{S}}[S_{\mathbf{G}}(n)]), & \text{on } A_n, \\ (\Gamma_{\mathbf{S}}[S_{\mathbf{G}}(n)] - \Gamma_{\mathbf{S}}[\lceil n/2 \rceil]) - (\Gamma_{\mathbf{R}}[\lfloor n/2 \rfloor] - \Gamma_{\mathbf{R}}[R_{\mathbf{G}}(n)]), & \text{on } A_n^c. \end{cases}$$

Finally, $G_n$ being the $(n - T_n)$th term of a bounded increment martingale is of order $\sqrt{n - T_n}$ ($O_p(n^{1/4})$) and, more importantly, normalized by $\sqrt{n - T_n}$ should converge to a normal limit. This is essentially the argument for Theorem 4.2, while Theorem 4.3 follows as an application of the above with a Borel–Cantelli type argument. Lemma 4.5 proves that $n^{-5/4}(M_{\mathbf{G}}(n) - M_{\mathbf{A}}(n)) \approx n^{-1/4}G_n$, Lemma 4.6 provides a uniform central limit theorem for the martingales behind $G_n$ and the proofs of the theorems complete the rest of the arguments.

### 4.3. *Proofs of the theorems.*

LEMMA 4.4. *For a greedy policy and the sequence of stopping times* $\{T_n\}_{n \geq 1}$ *defined in* (4.4), *we have*

$$(4.16) \quad \begin{aligned} \left(\frac{1}{n}\right) &[\tilde{N}_{\mathbf{R}}(R_{\mathbf{G}}(n)) - \tilde{N}_{\mathbf{R}}(R_{\mathbf{G}}(T_n))] \\ &\times [\tilde{N}_{\mathbf{S}}(S_{\mathbf{G}}(n)) - \tilde{N}_{\mathbf{S}}(S_{\mathbf{G}}(T_n))] \xrightarrow{\mathrm{d}} \left(\frac{\sigma_{\mathbf{R}}^2 + \sigma_{\mathbf{S}}^2}{8\mu}\right)\chi^2(1). \end{aligned}$$

*Moreover, we also have* $L^1$ *convergence.*

PROOF.    First, we will show that

$$(4.17) \quad \begin{aligned} &\left[\frac{\tilde{N}_{\mathbf{R}}(R_{\mathbf{G}}(n)) - \tilde{N}_{\mathbf{R}}(R_{\mathbf{G}}(T_n))}{R_{\mathbf{G}}(n) - R_{\mathbf{G}}(T_n)}\right] \\ &\times \left[\frac{\tilde{N}_{\mathbf{S}}(S_{\mathbf{G}}(n)) - \tilde{N}_{\mathbf{S}}(S_{\mathbf{G}}(T_n))}{S_{\mathbf{G}}(n) - S_{\mathbf{G}}(T_n)}\right] \xrightarrow{\mathrm{P}} \mu. \end{aligned}$$

Note that by inherent symmetry, Slutsky's theorem and the fact that both terms above are probability vectors, it suffices to show that

$$(4.18) \qquad \left[ \frac{N_{\mathbf{R}}[R_{\mathbf{G}}(n), j] - N_{\mathbf{R}}[R_{\mathbf{G}}(T_n), j]}{R_{\mathbf{G}}(n) - R_{\mathbf{G}}(T_n)} \right] \xrightarrow{\mathrm{P}} r_j, \qquad j \geq 1.$$

Now (4.13) combined with Theorem 5.1 of [6] on lag sums gives us (4.18) and, hence, (4.17).

Second, by Lemma 4.2 and Slutsky's theorem,

$$(4.19) \qquad \left( \frac{R_{\mathbf{G}}(n) - R_{\mathbf{G}}(T_n)}{n - T_n} \right) \left( \frac{S_{\mathbf{G}}(n) - S_{\mathbf{G}}(T_n)}{n - T_n} \right) \xrightarrow{\mathrm{P}} \frac{1}{4}.$$

Combining (4.17) and (4.19) with Lemma 4.1, and using Slutsky's theorem, we have (4.16). Now observe that (4.17) and (4.19) are nonnegative sequences bounded above by one. This together with the uniform integrability of the sequence $\{n^{-1}(n - T(n))^2\}_{n \geq 1}$ provided by Lemma 4.3 gives us the $L^1$ convergence. Hence, the proof. $\square$

PROOF OF THEOREM 4.1. Working on the set $\{R_{\mathbf{G}}(T_n) = \lceil n/2 \rceil\}$, we have

$$
\begin{aligned}
(4.20) \quad & M_{\mathbf{G}}(n) - M_{\mathbf{A}}(n) \\
& = (\tilde{N}_{\mathbf{R}}(R_{\mathbf{G}}(n)) - \tilde{N}_{\mathbf{R}}(\lceil n/2 \rceil)) \cdot \tilde{N}_{\mathbf{S}}(S_{\mathbf{G}}(T_n)) \\
& \quad + (\tilde{N}_{\mathbf{R}}(R_{\mathbf{G}}(n)) - \tilde{N}_{\mathbf{R}}(\lceil n/2 \rceil)) \\
& \qquad \times (\tilde{N}_{\mathbf{S}}(S_{\mathbf{G}}(n)) - \tilde{N}_{\mathbf{S}}(S_{\mathbf{G}}(T_n))) \\
& \quad - \tilde{N}_{\mathbf{R}}(\lceil n/2 \rceil) \cdot (\tilde{N}_{\mathbf{S}}(\lfloor n/2 \rfloor) - \tilde{N}_{\mathbf{S}}(S_{\mathbf{G}}(n))).
\end{aligned}
$$

The first term on the right-hand side can be written as

$$
\begin{aligned}
(4.21) \quad & (\tilde{N}_{\mathbf{R}}(R_{\mathbf{G}}(n)) - \tilde{N}_{\mathbf{R}}(\lceil n/2 \rceil) - [R_{\mathbf{G}}(n) - \lceil n/2 \rceil]\tilde{r}) \\
& \quad \times \tilde{N}_{\mathbf{S}}(S_{\mathbf{G}}(T_n)) + [R_{\mathbf{G}}(n) - \lceil n/2 \rceil]\Gamma_{\mathbf{S}}[S_{\mathbf{G}}(T_n)].
\end{aligned}
$$

The first expression in (4.21) has zero conditional expectation given $\mathcal{G}_{T_n}$ on the set $A_n$, as it is the $(n - T_n)$th term of a zero martingale. The argument for this assertion is similar to that found in the proof of Theorem 3.1. The third term on the right-hand side of (4.20) can be written as

$$
\begin{aligned}
(4.22) \quad & \tilde{N}_{\mathbf{R}}(\lceil n/2 \rceil) \cdot (\tilde{N}_{\mathbf{S}}(\lfloor n/2 \rfloor) - \tilde{N}_{\mathbf{S}}(S_{\mathbf{G}}(n)) - [R_{\mathbf{G}}(n) - \lceil n/2 \rceil]\tilde{s}) \\
& \quad + [R_{\mathbf{G}}(n) - \lceil n/2 \rceil]\Gamma_{\mathbf{R}}[\lceil n/2 \rceil].
\end{aligned}
$$

The first expression in (4.22) has a conditional expectation of zero on the set $A_n$ as it is independent of $\mathcal{G}_n (\supseteq \mathcal{G}_{T_n})$ and conditioned on $\mathcal{G}_n$, has zero mean. Using symmetry together with (4.21) and (4.22), we have

$$\frac{1}{n}|\mathbb{E}(M_\mathbf{G}(n) - M_\mathbf{A}(n))$$
$$- \mathbb{E}([\tilde{N}_\mathbf{R}(R_\mathbf{G}(n)) - \tilde{N}_\mathbf{R}(R_\mathbf{G}(T_n))] \cdot [\tilde{N}_\mathbf{S}(S_\mathbf{G}(n)) - \tilde{N}_\mathbf{S}(S_\mathbf{G}(T_n))])|$$
$$\leq \frac{1}{n}\mathbb{E}(|\Gamma_\mathbf{R}[R_\mathbf{G}(T_n)] - \Gamma_\mathbf{R}[S_\mathbf{G}(T_n)]||R_\mathbf{G}(n) - \lceil n/2\rceil|)$$
$$\leq \gamma\mathbb{E}\left(\frac{|R_\mathbf{G}(n) - \lceil n/2\rceil|}{n}\right) \to 0,$$

where the convergence to zero of the last term follows by Theorem 3.2 and the dominated convergence theorem. The theorem follows now by using Lemma 4.4. □

LEMMA 4.5.

$$\left(\frac{M_\mathbf{G}(n) - M_\mathbf{A}(n)}{n^{5/4}}\right) - \frac{G_n}{n^{1/4}} \xrightarrow{\text{a.s.}} 0 \qquad as\ n \to \infty.$$

PROOF. We start with a decomposition analogous to (4.20),

$$\frac{M_\mathbf{G}(n) - M_\mathbf{A}(n)}{n} - G_n$$
$$= \left(\frac{1}{n}\right)(\tilde{N}_\mathbf{R}(R_\mathbf{G}(n)) - \tilde{N}_\mathbf{R}(R_\mathbf{G}(T_n)))$$
$$\times (\tilde{N}_\mathbf{S}(S_\mathbf{G}(n)) - \tilde{N}_\mathbf{S}(S_\mathbf{G}(T_n)))$$
$$(4.23) \qquad + \mathrm{I}_{A_n}\left[(\tilde{N}_\mathbf{R}(R_\mathbf{G}(n)) - \tilde{N}_\mathbf{R}(R_\mathbf{G}(T_n))) \cdot \left(\frac{\tilde{N}_\mathbf{S}(S_\mathbf{G}(T_n))}{n} - \frac{\tilde{s}}{2}\right)\right.$$
$$\left. - (\tilde{N}_\mathbf{S}(\lfloor n/2\rfloor) - \tilde{N}_\mathbf{S}(S_\mathbf{G}(n))) \cdot \left(\frac{\tilde{N}_\mathbf{R}(R_\mathbf{G}(T_n))}{n} - \frac{\tilde{r}}{2}\right)\right]$$
$$+ \mathrm{I}_{A_n^c}\left[(\tilde{N}_\mathbf{S}(S_\mathbf{G}(n)) - \tilde{N}_\mathbf{S}(S_\mathbf{G}(T_n))) \cdot \left(\frac{\tilde{N}_\mathbf{R}(R_\mathbf{G}(T_n))}{n} - \frac{\tilde{r}}{2}\right)\right.$$
$$\left. - (\tilde{N}_\mathbf{R}(\lfloor n/2\rfloor) - \tilde{N}_\mathbf{R}(R_\mathbf{G}(n))) \cdot \left(\frac{\tilde{N}_\mathbf{S}(S_\mathbf{G}(T_n))}{n} - \frac{\tilde{s}}{2}\right)\right].$$

We now show that each term on the right-hand side of (4.23), upon division by $n^{1/4}$, converges almost surely to zero. For the first term, using

$$(\tilde{N}_\mathbf{R}(R_\mathbf{G}(n)) - \tilde{N}_\mathbf{R}(R_\mathbf{G}(T_n))) \cdot (\tilde{N}_\mathbf{S}(S_\mathbf{G}(n)) - \tilde{N}_\mathbf{S}(S_\mathbf{G}(T_n)))$$
$$\leq \frac{(n - T_n)^2}{4},$$

the result follows from Lemma 4.1. The second and third terms on the right-hand side of (4.23) are similar (by symmetry) and, hence, it suffices to deal solely with the second. We observe that similar arguments exist to show that each of the two expressions forming the second term, when divided by $n^{1/4}$, converges almost surely to zero. Hence, we give only the argument for the first. Since $R_{\mathbf{G}}(T_n) + S_{\mathbf{G}}(T_n) = T_n$, we observe that on $A_n$

$$(4.24) \quad \frac{n/2 - S_{\mathbf{G}}(T_n)}{n^{5/8}} = \underbrace{\left(\frac{n/2 - S_{\mathbf{G}}(T_n)}{\max(n - T_n, 1)}\right)}_{\text{bounded by 1}} \underbrace{\left[\frac{\max(n - T_n, 1)}{n^{5/8}}\right]}_{\xrightarrow{\text{a.s.}} 0 \text{ (Lemma 4.1)}} \xrightarrow{\text{a.s.}} 0.$$

Second, we define two sequences of random variables converging almost surely to zero. Let $\{U_n\}_{n \geq 1}$ be defined as

$$(4.25) \quad U_n := \underbrace{\left[\frac{R_{\mathbf{G}}(n) - R_{\mathbf{G}}(T_n)}{n - T_n}\right] \sqrt{\frac{S_{\mathbf{G}}(T_n)}{n}}}_{\text{bounded by 1}} \underbrace{\left(\frac{n - T_n}{n^{5/8}}\right)}_{\xrightarrow{\text{a.s.}} 0 \text{ (Lemma 4.1)}}$$

$$\times \underbrace{\sqrt{\frac{\log_2 S_{\mathbf{G}}(T_n)}{n^{1/4}}}}_{\xrightarrow{\text{a.s.}} 0} \qquad \forall n \geq 1$$

and $\{W_n\}_{n \geq 1}$ as

$$(4.26) \quad W_n := \underbrace{\left[\frac{R_{\mathbf{G}}(n) - R_{\mathbf{G}}(T_n)}{n - T_n}\right]}_{\text{bounded by 1}} \underbrace{\left(\frac{n - T_n}{n^{5/8}}\right)}_{\xrightarrow{\text{a.s.}} 0 \text{ (Lemma 4.1)}}$$

$$\times \underbrace{\left(\frac{n/2 - S_{\mathbf{G}}(T_n)}{n^{5/8}}\right)}_{\xrightarrow{\text{a.s.}} 0 \text{ by } (4.24)} \qquad \forall n \geq 1.$$

Third, using the above, we decompose the expression of interest as

$$\left(\frac{\tilde{N}_{\mathbf{R}}(R_{\mathbf{G}}(n)) - \tilde{N}_{\mathbf{R}}(R_{\mathbf{G}}(T_n))}{n^{1/4}}\right) \cdot \left(\frac{\tilde{N}_{\mathbf{S}}(S_{\mathbf{G}}(T_n))}{n} - \frac{\tilde{s}}{2}\right)$$

$$= U_n \left(\frac{\tilde{N}_{\mathbf{R}}(R_{\mathbf{G}}(n)) - \tilde{N}_{\mathbf{R}}(R_{\mathbf{G}}(T_n))}{R_{\mathbf{G}}(n) - R_{\mathbf{G}}(T_n)}\right)$$

$$\times \sqrt{\frac{S_{\mathbf{G}}(T_n)}{\log_2 S_{\mathbf{G}}(T_n)}} \left(\frac{\tilde{N}_{\mathbf{S}}(S_{\mathbf{G}}(T_n))}{S_{\mathbf{G}}(T_n)} - \tilde{s}\right)$$

$$- W_n \underbrace{\left( \frac{\tilde{N}_{\mathbf{R}}(R_{\mathbf{G}}(n)) - \tilde{N}_{\mathbf{R}}(R_{\mathbf{G}}(T_n))}{R_{\mathbf{G}}(n) - R_{\mathbf{G}}(T_n)} \right)}_{\text{bounded by } 1} \cdot \tilde{s}.$$

In view of (4.25) and (4.26), to show that the above converges almost surely to zero, it suffices to show that, with probability one,

$$\limsup_{n \to \infty} \left| \left( \frac{\tilde{N}_{\mathbf{R}}(R_{\mathbf{G}}(n)) - \tilde{N}_{\mathbf{R}}(R_{\mathbf{G}}(T_n))}{R_{\mathbf{G}}(n) - R_{\mathbf{G}}(T_n)} \right) \right.$$
$$\left. \times \sqrt{\frac{S_{\mathbf{G}}(T_n)}{\log_2 S_{\mathbf{G}}(T_n)}} \left( \frac{\tilde{N}_{\mathbf{S}}(S_{\mathbf{G}}(T_n))}{S_{\mathbf{G}}(T_n)} - \tilde{s} \right) \right| < \infty.$$

But this follows from the Cauchy–Schwarz inequality and Lemma A.1 since the first term is a probability vector. Hence, the proof. $\square$

The proof of Theorems 4.2 and 4.3 will require a uniform central limit theorem for a class of policies which can be described as *greedy with off-sets*. This is the content of the next lemma; below we describe some needed notation. Let $\mathbf{G}_\delta$, for $\delta \in [-\gamma, \gamma]$, be a policy satisfying

$$C_{\mathbf{G}_\delta}(n+1) = \begin{cases} 1, & \text{if } \Gamma_{\mathbf{S}}[S(n)] > \Gamma_{\mathbf{R}}[R(n)] + \delta, \\ 0, & \text{if } \Gamma_{\mathbf{S}}[S(n)] < \Gamma_{\mathbf{R}}[R(n)] + \delta, \end{cases} \qquad n = 1, 2, \ldots.$$

Let $\{X_{\mathbf{R}}^*(n)\}_{n \geq 1}$ and $\{X_{\mathbf{S}}^*(n)\}_{n \geq 1}$ denote two auxiliary sequences of i.i.d. random variables with $X_{\mathbf{R}}^* \overset{\mathrm{d}}{=} X_{\mathbf{R}}$ and $X_{\mathbf{S}}^* \overset{\mathrm{d}}{=} X_{\mathbf{S}}$, and let $\Gamma_{\mathbf{R}}^*(\cdot)$ and $\Gamma_{\mathbf{S}}^*(\cdot)$ denote their respective partial sums. For $\delta \in [-\gamma, \gamma]$, we define the sequence of random variables $\{Y_n^\delta\}_{n \geq 1}$ and $\{Z_n^\delta\}_{n \geq 1}$ as

$$Y_n^\delta := \left( \sqrt{\frac{2}{(\sigma_{\mathbf{R}}^2 + \sigma_{\mathbf{S}}^2)n}} \right) [\Gamma_{\mathbf{R}}[R_{\mathbf{G}_\delta}(n)] - \Gamma_{\mathbf{S}}^*(R_{\mathbf{G}_\delta}(n))], \qquad n \geq 1$$

and

$$Z_n^\delta := \left( \sqrt{\frac{2}{(\sigma_{\mathbf{R}}^2 + \sigma_{\mathbf{S}}^2)n}} \right) [\Gamma_{\mathbf{S}}[S_{\mathbf{G}_\delta}(n)] - \Gamma_{\mathbf{R}}^*(S_{\mathbf{G}_\delta}(n))], \qquad n \geq 1.$$

LEMMA 4.6.  *There exists a $K > 0$ such that*

$$\max_{\delta \in [-\gamma, \gamma]} \left[ \sup_{t \in \mathbb{R}} |\Pr(Y_n^\delta \leq t) - \Phi(t)|, \sup_{t \in \mathbb{R}} |\Pr(Z_n^\delta \leq t) - \Phi(t)| \right]$$
(4.27)
$$\leq K n^{-1/4} \log(n).$$

PROOF.  It suffices, by symmetry, to show that the first of the two expressions in (4.27) satisfies the bound. We use a filtration $\{\mathcal{H}_m\}_{m \geq 0}$ defined

for $m \geq 1$ as

$$\mathcal{H}_m = \mathcal{H}_0 \vee \sigma\langle L_{\mathbf{R}}(1), \ldots, L_{\mathbf{R}}(R_{\mathbf{G}_\delta}(m));$$
$$L_{\mathbf{S}}(1), \ldots, L_{\mathbf{S}}(S_{\mathbf{G}_\delta}(m)); X_{\mathbf{R}}^*(1), \ldots, X_{\mathbf{R}}^*(R_{\mathbf{G}_\delta}(m))\rangle,$$

with $\mathcal{H}_0$ containing all the information needed for randomization by $C_{\mathbf{G}_\delta}$. Also, we define, for a fixed $n \geq 1$,

$$D_m := \begin{cases} \dfrac{X_{\mathbf{R}}(R_{\mathbf{G}_\delta}(m)) - X_{\mathbf{S}}^*(R_{\mathbf{G}_\delta}(m))}{\sqrt{(\sigma_{\mathbf{R}}^2 + \sigma_{\mathbf{S}}^2)n/2}}, & \text{if } C_{\mathbf{G}_\delta}(m) = 1, \\ 0, & \text{if } C_{\mathbf{G}_\delta}(m) = 0, \end{cases} \quad m = 2, 3, \ldots, n,$$

with $D_1 := \frac{Y_1}{\sqrt{n}}$. By construction,

$$(4.28) \qquad \sum_{i=1}^{n} D_i = Y_n^\delta \quad \text{and} \quad \max_{i \leq n} D_i \leq n^{-1/2}\left(\frac{\gamma}{\sqrt{(\sigma_{\mathbf{R}}^2 + \sigma_{\mathbf{S}}^2)/2}}\right).$$

As $C_{\mathbf{G}_\delta}(m)$ is $\mathcal{H}_{m-1}$ measurable and both $X_{\mathbf{R}}(R_{\mathbf{G}_\delta}(m))$ and $X_{\mathbf{S}}^*(R_{\mathbf{G}_\delta}(m))$ are independent of $\mathcal{H}_{m-1}$, we have

$$\mathbb{E}(D_m|\mathcal{H}_{m-1}) = 0 \quad \text{and} \quad \mathbb{E}(D_m^2|\mathcal{H}_{m-1}) = \left(\frac{2}{n}\right)C_{\mathbf{G}_\delta}(m),$$

$$(4.29) \qquad\qquad\qquad\qquad\qquad\qquad\qquad\qquad\qquad 1 \leq m \leq n.$$

Hence, as $D_m^\delta$ is $\mathcal{H}_m$ measurable, $\{\sum_{i=1}^m D_i\}_{1 \leq m \leq n}$ is a martingale. As a consequence of (4.29), we have

$$V_n^2 := \sum_{i=1}^{n} \mathbb{E}(D_i^2|\mathcal{H}_{i-1}) = \left(\frac{2}{n}\right)R_{\mathbf{G}_\delta}(n).$$

This implies that

$$(4.30) \qquad \Pr(|V_n^2 - 1| > n^{-1/2}(\log(n))^2)$$
$$= \Pr\left(\left|\frac{R_{\mathbf{G}_\delta}(n) - n/2}{\sqrt{n}}\right| > \left(\frac{1}{2}\right)(\log(n))^2\right).$$

By an argument similar to that in the proof of Lemma 4.3, we get, analogous to (4.14), for $t \geq 1$ and $n \geq 4(1 + \frac{\gamma}{\mu})^2$,

$$(4.31) \qquad \Pr\left(\left[\frac{R_{\mathbf{G}_\delta}(n) - n/2}{\sqrt{n}}\right]^2 > t\right) \leq 2\exp\left\{-\left(\frac{\mu}{4\gamma}\right)^2\sqrt{t}\right\}.$$

Combining (4.30) and (4.31), we get

$$(4.32) \quad \Pr(|V_n^2 - 1| > n^{-1/2}(\log(n))^2) \leq \exp\left(\frac{32\gamma^2}{\mu^2}\right)\left(\frac{1}{n}\right) \qquad \forall n \geq 1.$$

Inequalities in (4.28) and (4.32) imply that the two conditions of Lemma A.2 are satisfied in our case. Hence, we have (4.27), for some $K$ free of $\delta$. Hence, the proof.   □

PROOF OF THEOREM 4.2.  In view of Lemma 4.5, it suffices to derive the weak limit of $n^{-1/4}G_n$. We start by observing that, for $u \in \mathbb{R}$ on $A_n$,

$$\Pr\left(\frac{(\Gamma_{\mathbf{R}}[R_{\mathbf{G}}(n)] - \Gamma_{\mathbf{R}}[\lceil n/2 \rceil]) - (\Gamma_{\mathbf{S}}[\lfloor n/2 \rfloor] - \Gamma_{\mathbf{S}}[S_{\mathbf{G}}(n)])}{\sqrt{(\sigma_{\mathbf{R}}^2 + \sigma_{\mathbf{S}}^2)(n - T_n)/2}} \le u \,\Big|\, \mathcal{G}_{T_n}\right)$$

$$= \Pr(Y_{n-T_n}^{\Delta_n} \le u)$$

and on $A_n^c$,

$$\Pr\left(\frac{(\Gamma_{\mathbf{S}}[S_{\mathbf{G}}(n)] - \Gamma_{\mathbf{S}}[\lceil n/2 \rceil]) - (\Gamma_{\mathbf{R}}[\lfloor n/2 \rfloor] - \Gamma_{\mathbf{R}}[R_{\mathbf{G}}(n)])}{\sqrt{(\sigma_{\mathbf{R}}^2 + \sigma_{\mathbf{S}}^2)(n - T_n)/2}} \le u \,\Big|\, \mathcal{G}_{T_n}\right)$$

$$= \Pr(Z_{n-T_n}^{\Delta_n} \le u),$$

where $\Delta_n := \Gamma_{\mathbf{R}}[R_{\mathbf{G}}(T_n)] - \Gamma_{\mathbf{S}}[S_{\mathbf{G}}(T_n)]$. This, along with Lemma 4.6, leads to

$$\left|\Pr\left(\frac{G_n}{\sqrt{0.125(\sigma_{\mathbf{R}}^2 + \sigma_{\mathbf{S}}^2)(n - T_n)}} \le u\right) - \Phi(u)\right|$$

$$= \left|\int_{A_n} \Pr(Y_{n-T_n}^{\Delta_n} \le u)\, dP + \int_{A_n^c} \Pr(Z_{n-T_n}^{\Delta_n} \le u)\, dP - \Phi(u)\right|$$

$$\le K\mathbb{E}[\min[1, (n - T_n)^{-1/4} \log(n - T_n)]] \to 0 \qquad \text{as } n \to \infty.$$

In other words, we have shown that

$$\frac{G_n}{\sqrt{n - T_n}} \xrightarrow{\mathrm{d}} \mathrm{N}\left(0, \frac{\sigma_{\mathbf{R}}^2 + \sigma_{\mathbf{S}}^2}{8}\right) \qquad \text{as } n \to \infty.$$

This with the asymptotic independence between the terms on the right of

$$\left(\frac{G_n}{n^{1/4}}\right) = \left(\frac{G_n}{\sqrt{n - T_n}}\right)\left(\sqrt{\frac{n - T_n}{\sqrt{n}}}\right)$$

and Lemma 4.1 completes the proof.   □

PROOF OF THEOREM 4.3.  In view of Lemma 4.5, to prove (4.2), it suffices to show that

$$\liminf \frac{G_n}{(n \log_2(n))^{1/4}} = -\infty \quad \text{and} \quad \limsup \frac{G_n}{(n \log_2(n))^{1/4}} = \infty.$$

Due to the similarity of the arguments, we prove only the latter. Toward this end, we define a sequence of stopping times $\{T_n^*\}_{n\geq 1}$ as

$$T_n^* = \begin{cases} \inf\{T_{2k} \geq T_{n-1}^* | k \geq 2; R_{\mathbf{G}}(T_{2k}) = k; \\ \qquad 2k - T_{2k} \geq \sigma_{R_{\mathbf{G}}}\sqrt{2k\log_2(2k)}\}, & n \text{ odd}, \\ \inf\{k \geq T_{n-1}^* | R_{\mathbf{G}}(k) = S_{\mathbf{G}}(k)\}, & n \text{ even}. \end{cases}$$

The stopping times are easily checked to be well defined using the definition of $\{T_n\}_{n\geq 1}$ and Theorem 3.3. Let $C > 0$ be an arbitrary constant and $\{B_i\}_{i\geq 1}$ be a sequence of events defined by

$$B_i := \left\{ \frac{G_{2R_{\mathbf{G}}(T_{2i-1}^*)}}{(2\sigma_{R_{\mathbf{G}}}^2 R_{\mathbf{G}}(T_{2i-1}^*)\log_2(2R_{\mathbf{G}}(T_{2i-1}^*)))^{1/4}} > C \right\}, \qquad i = 1, 2, \dots.$$

Also let $\{\mathcal{H}_i\}_{i\geq 0}$ be a filtration with $\mathcal{H}_i := \mathcal{G}_{T_{2i+1}^*}$ for $i \geq 0$. By construction, $B_i \in \mathcal{H}_i$ for $i \geq 1$ and, moreover, applying Lemma 4.6 as in Theorem 4.2, we have

$$\Pr(B_i|\mathcal{H}_{i-1}) \geq \Pr\left( \frac{G_{2R_{\mathbf{G}}(T_{2i-1}^*)}}{\sqrt{2R_{\mathbf{G}}(T_{2i-1}^*) - T_{2i-1}^*}} > C \Big| \mathcal{H}_{i-1} \right)$$

$$> \frac{1 - \Phi(C)}{2} > 0, \qquad \text{for large } i.$$

Now Lemma A.3 implies that, with probability one, $B_i$ occurs infinitely often. This completes the proof of (4.2).

To show (4.3), it suffices to look at the subsequence of even epochs. If $R_{\mathbf{G}}(2n) = n$, then $M_{\mathbf{G}}(2n) = M_{\mathbf{A}}(2n)$. Suppose that $R_{\mathbf{G}}(2n) = n + K_n > n$ (swap $\mathbf{R}$ with $\mathbf{S}$ in the contrary). This leads to

$$\begin{aligned} M_{\mathbf{G}}(2n) - M_{\mathbf{A}}(2n) \leq{} & n(\Gamma_{\mathbf{R}}[n + K_n] - \Gamma_{\mathbf{R}}[n]) \\ & + (\tilde{N}_{\mathbf{S}}(n) - n\tilde{s}) \cdot (\tilde{N}_{\mathbf{R}}(n + K_n) - \tilde{N}_{\mathbf{R}}(n)) \\ & - n(\Gamma_{\mathbf{S}}[n] - \Gamma_{\mathbf{S}}[n - K_n]) \\ & - (\tilde{N}_{\mathbf{R}}(n) - n\tilde{r}) \cdot (\tilde{N}_{\mathbf{S}}(n) - \tilde{N}_{\mathbf{S}}(n - K_n)). \end{aligned}$$

All statements which follow are to be understood as holding eventually, with probability one. By Theorem 3.3, $K_n \leq B(n\log_2 n)^{1/2}$ for some constant $B$ and by Lemma A.1, the components of $(\tilde{N}_{\mathbf{S}}(n) - n\tilde{s})$ are uniformly bounded by $n^{9/16}$. Since the components of $(\tilde{N}_{\mathbf{R}}(n + K_n) - \tilde{N}_{\mathbf{R}}(n))$ are bounded by $K_n$, we have the second term above is of order at most $n^{9/8}$. The same holds true for the fourth term. By Theorem 5.1 of [6] on lag sums, we have, for $0 \leq k \leq B(n\log_2 n)^{1/2}$ and for some constant $C$,

$$\begin{aligned} n(\Gamma_{\mathbf{R}}[n + k] - \Gamma_{\mathbf{R}}[n]) &- n(\Gamma_{\mathbf{S}}[n] - \Gamma_{\mathbf{S}}[n - k]) \\ &\leq n(k\mu + C\sqrt{k\log n}) - n(k\mu - C\sqrt{k\log n}) \\ &\leq 2Cn^{5/4}\sqrt{\log n}(\log_2 n)^{1/4}. \end{aligned}$$

This implies that $M_\mathbf{G}(2n) - M_\mathbf{A}(2n) \leq 3Cn^{5/4}\sqrt{\log n}(\log_2 n)^{1/4}$. Similarly, the difference $M_\mathbf{G}(2n) - M_\mathbf{A}(2n)$ can be bounded from below. Hence, the proof. $\quad\square$

## APPENDIX

The first lemma, on the rate of $l^2$ convergence of empirical probabilities to the true probabilities, derives from [14] and is included here for the reader's convenience.

LEMMA A.1.  *Let* $\{Z_i\}_{i\geq 1}$ *be a random sample from a discrete distribution described by*

$$p_j := \Pr(Z_1 = z_j), \qquad j = 1, 2, \dots \text{ with } \sum_{j\geq 1} p_j = 1.$$

*Defining the empirical probability vector* $(p_j^n)_{j\geq 1}$, *for* $n \geq 1$, *by*

$$p_j^n = \left(\frac{1}{n}\right)\#\{1 \leq k \leq n : X_k = z_j\}, \qquad j = 1, 2, \dots,$$

*we have*

$$\limsup\left(\frac{n}{\log_2(n)}\right)\sum_{j=1}^{\infty}(p_j^n - p_j)^2 \xrightarrow{\text{a.s.}} C < \infty.$$

PROOF.  Without loss of generality, we assume that $z_j = j$, for $j \geq 1$. Defining the kernel $h(\cdot, \cdot)$ by

$$h(i,j) = I_{\{i=j\}} - (p_i + p_j) + \sum_{k\geq 1}(p_k)^2, \qquad i, j = 1, 2, \dots,$$

we have

$$
\begin{aligned}
\left(\frac{n}{\log_2(n)}\right)\sum_{j=1}^{\infty}(p_j^n - p_j)^2 &= \left(\frac{2}{n\log_2(n)}\right)\sum_{1\leq i<j\leq n} h(X_i, X_j) \\
&\quad + \left(\frac{1}{\log_2(n)}\right)\left(1 + \sum_{k\geq 1}(p_k)^2\right) \\
&\quad - \left(\frac{2}{\log_2(n)}\right)\left(\frac{\sum_{i=1}^{n} p_{X_i}}{n}\right).
\end{aligned}
$$

It is easy to check that the first term on the right is a canonical $U$-statistic of order 2. By [2], we have

$$\limsup\left(\frac{2}{n\log_2(n)}\right)\sum_{1\leq i<j\leq n} h(X_i, X_j) \xrightarrow{\text{a.s.}} C < \infty$$

for some constant $C$. The third term on the right converges to zero in the almost sure sense by the usual SLLN. Hence, the proof. $\quad\square$

The following uniform central limit theorem for martingales is a restatement of Theorem 3.7 of [5] for a noncanonical filtration [see Remark (ii) on page 84 following the theorem].

LEMMA A.2. *Let* $\{S_i = \sum_1^i X_j, \mathcal{H}_i, 1 \leq i \leq n\}$ *be a zero-mean martingale. Let*

$$V_i^2 = \sum_1^i \mathbb{E}(X_j^2 | \mathcal{H}_{j-1}), \qquad 1 \leq i \leq n$$

*and suppose that*

$$\max_{i \leq n} |X_i| \leq n^{-1/2} M \qquad a.s.$$

*and*

$$\Pr(|V_n^2 - 1| > 9M^2 D n^{-1/2} (\log n)^2) \leq C n^{-1/4} \log n$$

*for constants* $M$, $C$ *and* $D(\geq e)$. *Then for* $n \geq 2$,

$$\sup_{x \in \mathbb{R}} |\Pr(S_n \leq x) - \Phi(x)| \leq K n^{-1/4} \log n,$$

*where* $K$ *is a universal function of* $M$, $C$ *and* $D$.

The last result is a conditional Borel–Cantelli lemma which appears as Theorem 2.8.5 in [17].

LEMMA A.3. *Let* $\{B_i, i \geq 1\}$ *be a sequence of events and* $\{\mathcal{H}_i, i \geq 1\}$ *an increasing sequence of* $\sigma$-*fields such that* $B_i \in \mathcal{H}_i$ *for each* $i \geq 1$. *Then*

$$\{B_i \ i.o.\} = \left\{ \sum_{i=1}^{\infty} \Pr(B_i | \mathcal{H}_{i-1}) = \infty \right\},$$

*that is,* $\sum_{i=1}^{\infty} \Pr(B_i | \mathcal{H}_{i-1}) < \infty$ *implies the* $B_i$ *occur at most finitely often and* $\sum_{i=1}^{\infty} \Pr(B_i | \mathcal{H}_{i-1}) = \infty$ *implies the* $B_i$ *occur infinitely often.*

**Acknowledgments.** We thank Professor Richard Dykstra for fruitful discussions and Professor Ramon Lawrence for bringing the problem to our attention. Also, we thank an associate editor and a referee whose comments have led to a significant enhancement of our results.

DEPARTMENT OF STATISTICS AND ACTUARIAL SCIENCE
UNIVERSITY OF IOWA
IOWA CITY, IOWA 52242
USA
E-MAIL: rrusso@stat.uiowa.edu
        nshyamal@stat.uiowa.edu